\documentclass[12pt]{amsart}
\usepackage{hyperref}\hypersetup{colorlinks=true, linkcolor=purple, filecolor=magenta, urlcolor=cyan}
\usepackage{amsfonts}
\usepackage{dsfont}
\usepackage{mathrsfs} 
\usepackage{amssymb, amsthm, changes}
\usepackage{amsxtra}
\usepackage{amstext}
\usepackage[english]{babel}
\AtBeginDocument{ 
   \def\MR#1{}
}

\usepackage{enumitem}
\newtheorem{theorem}{Theorem}[section]
\newtheorem{lemma}[theorem]{Lemma}

\newtheorem{proposition}[theorem]{Proposition}

\newtheorem{remark}[theorem]{Remark}

\newtheorem{definition}{Definition}

\newtheorem{assumption}{Assumption}

\newcommand{\ZZ}{\mathbf{Z}}
\newcommand{\E}{\mathbb{E}}
\renewcommand{\P}{\mathbb{P}}

\newcommand{\supp}{\mathrm{supp}\,}

\newcommand{\X}{\mathbf{X }}

\newcommand{\R}{\mathbb{R}}
\newcommand{\Z}{\mathbb{Z}}
\newcommand{\N}{\mathbb{N}}
\def\eps{\varepsilon}
\def\phi{\varphi}

\def\l{\ell}
\def\1{\mathds{1}}
\def\1vec{\mathbf{1}}
\def\0{\mathbf{0}}

\usepackage{tcolorbox}

\begin{document}

\title{Speed of excited random walks with long backward steps}

\author{Tuan-Minh Nguyen}
\address{School of Mathematical Sciences, Monash University, Victoria 3800, Australia}
\email{tuanminh.nguyen@monash.edu}

\begin{abstract}We study a model of multi-excited random walk with non-nearest neighbour steps on $\Z$, in which the walk can jump from a vertex $x$ to either $x+1$ or $x-i$ with $i\in \{1,2,\dots,L\}$, $L\ge 1$. We first point out the multi-type branching structure of this random walk and then prove a limit theorem for a related multi-type Galton-Watson process with emigration, which is of independent interest. Combining this result and the method introduced  by Basdevant and Singh [Probab. Theory Related Fields  (2008), {\bf 141}(3-4)], we extend their result (w.r.t. the case $L=1$) to our model. More specifically, we show that in the  regime of transience to the right, the walk has positive speed if and only if the expected total drift $\delta>2$. This confirms a special case of a conjecture proposed by Davis and Peterson.
\end {abstract}



\maketitle
\tableofcontents

\section{Introduction}

Excited random walk is a model of non-markovian random walk in a cookie environment, in which the walker consumes a cookie (if available) upon reaching a site and makes a jump with transition law  dynamically depending on the number of remaining cookies at its current position. The model of nearest-neighbour excited random walks has been extensively studied in recent years. Benjamini and Wilson \cite{BW2003} first studied once-excited random walks with a focus on higher-dimensional integer lattice. Later, Zerner \cite{Z2005}  extended this model to multi-excited random walks and established a  criterion for recurrence/transience of the model on $\Z$.   There are also notable results for asymptotic behaviour of the multi-excited model including criteria non-ballisticity/ballisticity~\cite{DS2008} as well as characterization of the limit distribution in such specific regimes ~\cite{DK2012, KM2011, Kosygina2014}. See also \cite{Hofstad2010, Menshikov2012, KP2017, KMP2022}. For a literature review, we refer the reader to \cite{KZ2013}.

\subsection{Description of the model and the main result}
We define a non-nearest-neighbour random walk $\X := (X_n)_{n\ge0}$, which describes the position of a particle moving in a cookie environment on the integers $\Z$ as follows. For any integer $n$, set $[n] = \{1, 2, \ldots, n\}$.
Let $M$ and $L$  be positive integers, $\nu$ and $(q_j)_{j\in [M]}$ be probability measures on $\Lambda:=\{-L, -L+1,\dots,-1,1 \}$.  Initially, each vertex in $\Z$ is assigned  a stack of $M$ cookies and we set $X_0=0$.  Suppose that $X_n =x$ and by time $n$ there are exactly  remaining $M-j+1$ cookie(s) at site $x$ with some $j\in [M]$. Before the particle  jumps to a different location, it eats one cookie and jumps to site $x+i$, $i\in \Lambda$, with probability $q_j(i)$. On the other hand, if the stack of cookies at $x$ is empty then it jumps to site $x+i$, $i\in\Lambda$ with probability $\nu(i)$. More formally, denote by $(\mathcal{F}_n)_{n\ge 0}$ the natural filtration of $\X$.  For each $i\in \Lambda$
$$\P\left(X_{n+1}=X_n+i | \mathcal{F}_n \right)= \omega(\mathcal{L}(X_n,n),i)$$
where $\mathcal{L}(x,n)=\sum_{k=0}^n \mathds{1}_{\{X_k=x\}}$ is the number of visits to vertex $x\in\mathbb{Z}$ up to time $n$, and $\omega: \N\times \Lambda \to [0,1]$ is the \textit{cookie environment} given by
$$\omega(j,i) =\left\{\begin{array}{ll} q_{j}(i), \quad & \text{if}\ 1\le j \le M,\\
\nu(i), \quad  & \text{if}\ j >  M. \end{array}\right.$$

Throughout this paper, we make the following assumption. \vspace{5pt}

\begin{tcolorbox}
\begin{assumption}\label{Asum:A}
{\begin{itemize}\item[i.] The distribution $\nu$ has zero mean. 
\item[ii.]  For each $j\in[M]$, the distribution $q_j$ has positive mean and $q_j(1)<1$.
\end{itemize}}
\end{assumption}
\end{tcolorbox}
We call the process $\X$ described above $(L,1)$ non-nearest neighbors excited random walk ($(L,1)$-ERW, for brevity). 
It is worth mentioning that $(L,1)$-ERW is a special case of excited random walks with non-nearest neighbour steps considered by Davis and Peterson in \cite{DP2017} in which the particle can also jump to  non-nearest neighbours on the right  and $\Lambda$ can be an unbounded subset of $\Z$. In particular, Theorem 1.6 in \cite{DP2017} implies that the  process studied in this paper is 
\begin{itemize}
\item  transient to the right if the \textit{expected total drift} $\delta$, defined as
\begin{equation}\label{eq:totdr} \delta:=\sum_{j=1}^M\sum_{\l\in \Lambda} \l q_j(\l) 
\end{equation}
 is larger than 1, and 
 \item recurrent if $\delta\in [0,1]$. 
 \end{itemize}
Additionally, Davis and Peterson conjectured that the limiting speed of the random walk exists if $\delta>1$ and it is positive when $\delta>2$. (see Conjecture 1.8 in \cite{DP2017}). 

Recently, a sufficient condition for once-excited random walks with long forward jumps to have positive speed has been shown in \cite{CHN2021}. However, the coupling method introduced in \cite{CHN2021} seems to not be applicable to models of multi-excited random walks. 


In the present paper, we verify Davis-Peterson conjecture for $(L,1)$-ERW. More precisely, we show that $\delta>2$ is  a sufficient and necessary condition for $(L,1)$-ERW to have  positive limiting speed, under Assumption~\ref{Asum:A}.
\begin{tcolorbox}
\begin{theorem}\label{th:main} Under Assumption~\ref{Asum:A},
\begin{itemize}
\item[(a)] if  $\delta >1$ the speed of $(L,1)$-ERW  $\X$ exists, i.e. 
$X_n/n$ converges a.s. to a non-negative constant $v$, and
\item[(b)]  if $\delta >2$ we have that $v >0$. If $\delta \in (1, 2]$ then $v =0$.
\end{itemize}
\end{theorem} 
\end{tcolorbox}
\subsection{Summary of the proof of Theorem \ref{th:main}.}
Our proof strategy relies on the connection between non-nearest neighbor excited random walks and multi-type branching processes with migration.  The idea can be traced back to the branching structures of nearest-neighbor excited random walks \cite{DS2008} and random walks in a random environment (see e.g. \cite{Kesten} and \cite{HW2016}). In the present paper, we introduce  multi-type branching process with emigration and develop techniques from \cite{DS2008} to deal with various higher dimensional issues in our model.

The remaining parts of the paper are organized as follows. We first describe in Section \ref{sec:branching} the multi-type branching structure of the number of backward jumps. This branching structure is formulated by a multi-type branching process with (random) migration $\ZZ$ defined in Proposition \ref{Markovian}. In Section~\ref{sec:emigtation}, we next demonstrate a limiting theorem for a class of critical multi-type Galton-Watson processes with emigration. We believe that this result is of independent interest. In section \ref{sec:speed}, we derive  a functional equation related to limiting distribution of $\ZZ$ (Propositions \ref{prop:1} and \ref{prop:2}). Combining these results together with a coupling between $\ZZ$ and a critical multi-type branching process with emigration (which is studied in Section~\ref{sec:emigtation}), we deduce the claim of Theorem \ref{th:main}.

It is worth mentioning that the techniques introduced in this paper is unfortunately not applicable to the case of excited random walks having {non-nearest-neighbour} jumps to the right. We refer the reader to \cite{CHN2021} for a recent work studying the speed of once-excited random walks with long forward steps.   

\section{Multi-branching structure of excited random walks}\label{sec:branching}
{For any pair of functions $f$ and $g$ of one real or discrete variable, we write
\begin{itemize}
\item $f(x) \sim g(x)$ as $x\to x_0$ if
$\lim_{x \to x_0} {f(x)}/{g(x)} =1,$
\item $f(x)=O(g(x))$ as $x\to x_0$ if $\limsup_{x\to x_0}|f(x)/g(x)| <\infty $ and
\item  $f(x)=o(g(x))$ as $x\to x_0$ if $\lim_{x\to x_0}f(x)/g(x)=0$.
\end{itemize}}

Denote  $\Z_+:=\{0,1,2,\dots\}$ and $\N=\Z_+\setminus\{0\}$. For any $m,n \in \Z$, $m\le n$, set $[m,n]_\Z:= \{m, \ldots, n\}$ and $[n]=\{1,2,\dots,n\}$.

For each $n\in \Z_+$, let $T_n=\inf\{ k\ge0: X_k=n\}$ be the first hitting time of site $n$. For $i\le n-1$, define $V_i^n=(V_{i,1}^n, V_{i,2}^n,\dots, V_{i,L}^n)$ where for $\l\in [L]$,
$$V_{i,\l}^n=\sum_{k=0}^{T_n-1}\mathds{1}_{\{X_k>i,X_{k+1}=i-\l+1 \}}$$
stands for the number of backward jumps from a site in  the set $ i + \N$  to site $i-\l+1$ before time $T_n$. Notice that $T_n$ is equal to the total number of forward and backward jumps before time $T_n$. In particular, the number backward jumps to site $i$ before time $T_n$ is equal to $V_{i,1}^n$. On the other hand, {between} two consecutive forward jumps from $i$ to $i+1$, there is exactly one backward jump from $ i + \N$ to $ i - \Z_+$. Furthermore, for $0\le  i\le n-1$, before the first backward jump from $ i + \N$ to $ i - \Z_+$, the walk must have its first forward jump from $i$ to $i+1$. Therefore the number of forward jumps from $i$ to $i+1$ before time $T_n$ is equal to 
$1_{\{0\le i\le n-1\}}+\sum_{\l=1}^LV_{i,\l}^n$. As a result, we obtain
$$ T_n= n +2\sum_{ -\infty<i\le n-1}  V_{ i,1}^n+\sum_{-\infty<i\le n-1} \sum_{\l=2}^L V_{i,\l}^n.$$
Assume from now that $(X_n)_{n\ge 0}$ is transient to the right. Notice that the walk spends only a finite amount of time on $-\N$ and thus
\begin{equation}\label{hitting}
T_n\sim  n +2\sum_{i= 0}^{n-1} V_{i,1}^{n}+\sum_{i= 0}^{n-1}\sum_{\l=2}^L V_{i,\l}^n \quad \text{as\ } n \to\infty, \qquad \mbox{a.s..}
\end{equation}
It is worth mentioning that the above hitting time decomposition was mentioned by Hong and Wang \cite{HW2016}, in which they studied random walks in random environment with non-nearest-neighbour jumps to the left. The idea can be traced back to the well-known Kesten-Kozlov-Spitzer hitting time decomposition for nearest-neighbour random walks in random environments \cite{Kesten}.  

Let $(\xi_n)_{n\in\N}$ be a sequence of independent random unit vectors such that
the distribution of $\xi_n=(\xi_{n,1},\xi_{n,2},\dots,\xi_{n,L+1})$ is given by $$\P(\xi_n=\mathbf{e}_{\ell})=\left\lbrace \begin{matrix}& q_{n}(-\ell),\quad &\text{if }1\le n \le M \text{ and } 1\le \ell\le L,  \\
& q_n(1),\quad & \text{if }1\le n \le M \text{ and } \ell=L+1,\\
&\nu(-\ell),\quad & \text{if } n>M,1\le \ell\le L,\\
& \nu(1) ,\quad & \text{if } n>M,  \ell= L+1  \end{matrix}\right.$$
where $\mathbf{e}_{\ell}$ with $\ell \in [L+1]$ is the standard basis of $\R^{L+1}.$ 
 If $\xi_n=\mathbf{e}_{\l}$ with $\l\in [L]$, we say that the outcome of the $n$-th experiment is an $\ell$-th type failure. Otherwise, if $\xi_n=\mathbf e_{L+1}$, we say that it is a success.   

For $m\in \Z_+$, we define the random vector $A(m)=(A_1(m),\dots, A_L(m))$  such that for $\l\in [L]$,
\begin{equation}\label{def:Aj}A_{\l}(m):=\sum_{i=1}^{\gamma_m}\xi_{i,\l}, \quad \text{with} \quad \gamma_m=\inf\Big\{n\ge 1: \sum_{i=1}^n \xi_{i,L+1}=m+1\Big\}.
\end{equation}
In other words, the random variable $A_{\l}(m)$ is the total number of $\ell$-th type failures before obtaining $m+1$ successes. 

{Let $(A^{(n)}(m))_{m\in \Z_+}$, with $n\in \N$, be i.i.d. copies of the process $(A(m))_{m\in \Z_+}$. We define a $L$-dimensional process $\ZZ =(Z_n)_{n\ge0}=(Z_{n,1},Z_{n,2},\dots,Z_{n,L})_{n\ge0}$ such that $Z_0\in \Z_{+}^L$ is independent of $(A^{(n)}(m))_{n\in\N, m\in \Z_+}$ and for $n\ge 1$, \begin{align}\label{def.Z}Z_{n}=A^{(n)}\left(|Z_{n-1}|\right)+ (Z_{n-1,2},Z_{n-1,3},\dots,Z_{n-1,L},0)\end{align}
where for $z=(z_1, z_1,\dots,z_{L})\in \R^{L}$, we denote $|z|=z_1+z_2+\dots+z_{L}$. Therefore $\ZZ=(Z_n)_{n\ge 0}$ is a Markov chain in $\Z_{+}^L$ and its transition law is given by}
\begin{align*}&\P\left(Z_{n+1}=(k_{1},k_{2},\dots,k_L)  \ \big| \  Z_{n}=(j_1,j_2,\dots, j_L)\right)\\
\nonumber & = \P\left( A_1\Big(\sum_{\l=1}^L j_{\l}\Big)=k_1-j_2,\dots, A_{L-1}\Big(\sum_{\l=1}^L j_{\l}\Big)=k_{L-1}-j_L,A_L\Big(\sum_{\l=1}^L j_{\l}\Big)=k_{L}  \right)
\end{align*}
for $(k_{1},k_{2},\dots,k_L), (j_{1},j_{2},\dots,j_L)\in \Z_+^L$.
\begin{proposition}\label{Markovian} {Assume  that $Z_{0,\l}=0$ for all $\l \in [L]$. Then for each $n\in \N $, we have that 
$(V_{n-1}^n,V_{n-2}^n,\dots, V_{0}^n)$ has the same distribution as $(Z_0,Z_1,\dots,Z_{n-1})$.} 
\end{proposition}

\begin{proof}

A backward jump is called $\l$-th type  of level $i$ if it is a backward jump from a site in $i + \N$  to site $i-\l+1$. Recall that {$V_{i,\l}^n$} is the number of $\l$-th type backward jumps of level $i$ before time  $T_n$. Assume that $\{V_{i}^n=(V_{i,1}^n,\dots,V_{i,L}^n)=(j_1,j_2,\dots,j_L)\}$. The number of forward jumps from $i$ to $i+1$ before time  $T_n$ is thus equals to $1+\sum_{\l=1}^L V_{i,\l}^n=1+\sum_{\l =1}^L j_\ell$. 

For each $i\in \Z$, denote by $T_{i}^{(k)}$ the time for $k$-th forward jump from $i-1$ to $i$ and also set $T^{(0)}_i=0$. We have that  $T_{i}^{(1)} = T_i$. Moreover, as the process $\X$ is transient, we have that only finitely many  $(T_i^{(k)})_k$  are finite,  and conditioning on $\{V_{i}^n=(j_1,j_2,\dots,j_L)\}$, we have that $T_i^{k} < \infty$ for $k \le 
1+\sum_{s=1}^L j_s$. 

Note that $V_{i-1,\l}^n$, i.e. the number of $\l$-th type backward jumps of level $i-1$ before time $T_n$, is equal to the sum of number of $\l$-th type backward jumps of level $i-1$ during $[T_{i+1}^{(k-1)},T_{i+1}^{(k)}-1]_{\Z}$ for $k \in [1+\sum_{k=1}^L j_k]$. 
 
 By the definition of $T_{i}^{(k)}$,  the walk will   visit $i$ at least once during the time interval $[T_{i+1}^{(k-1)},T_{i+1}^{(k)}-1]_{\Z}$. Whenever the walk visits $i$, it will make a forward jump from $i$ to $i+1$ (which corresponds to a success) or a backward jump from $i$ to $i-\ell$, i.e. a $\l$-th type jump of level $i-1$, with $\l\in [L]$. If the latter happens, then $i$ will be visited again during $[T_{i+1}^{(k-1)},T_{i+1}^{(k)}-1]_{\Z}$. Moreover, an $\l$-th type backward jump of level $i$ is also an $(\l-1)$-th type backward jump of level $i-1$. Thus conditionally on $\{(V_{i,1}^n,\dots,V_{i,L}^n)=(j_1,j_2,\dots,j_L)\}$,
the random vector $(V_{i-1,1}^n,V_{i-1,2}^n,\dots, V_{i-1,L}^n)$ has the same distribution as 
$$\left( A_1\Big(\sum_{\l=1}^L j_{\l}\Big)+j_2,\dots , A_{L-1}\Big(\sum_{\l=1}^L j_{\l}\Big)+j_L ,  A_L\Big(\sum_{\l=1}^L j_{\l}\Big)\right).$$
\end{proof}
{Recall that $M$ is the total number of cookies initially placed on each site.} By the definition of sequence $(A(m))_{m\in\Z_+}$ given in \eqref{def:Aj}, we can easily obtain the following. 
\begin{proposition}\label{iden.A} For $m\ge M-1$, we have
\begin{align*} 
A({m})=A({M-1})+ \sum_{k=1}^{m-M+1} \eta_{k}
\end{align*}
where $(\eta_k)_{k\ge1}$ are i.i.d. random vectors independent of $A({M-1})$ with multivariate geometrical law 
\begin{align}\label{geom}
\P\left(\eta_1=(i_1,i_2,\dots,i_{L})\right)= \frac{\nu(1) }{(i_1+i_2+\dots+i_L)!} \prod_{k \in [L]}  i_k ! \nu(-k)^{i_{k}}.
\end{align}
In the above formula, we use the convention that $0^0=1$.
\end{proposition}
\begin{remark}\label{rem:branching}
The multivariate Markov chain $\mathbf Z$ defined in \eqref{def.Z} can be interpreted as a multi-type branching process with (random) migration as follows. Let $M'$ be a fixed integer such that $M'\ge M-1$ and suppose that $Z_{n-1}=j=(j_1,j_2,\dots,j_L)$. Then
\begin{itemize}[leftmargin=15pt]
 \item If $|j|=j_1+j_2+\dots +j_L\ge M'$, we have {$$Z_{n}=A^{(n)}({M'})+\sum_{k=M'-M+2}^{|j|-M+1}\eta_k^{(n)}+ \widetilde{j}$$ where $\widetilde{j}:=(j_2,\dots, j_{L},0)$; $\eta_k^{(n)}$ with $k\in\N$ and $n\in\N$ are i.i.d. random vectors with the multivariate geometrical law defined in \eqref{geom}; and for each $n\in\N$, $(\eta_k^{(n)})_{k\in \N}$ is independent of $A^{(n)}(M')$ and $(Z_{k})_{0\le k\le n-1}$.} In this case, there is an emigration of $M'$ particles (each particle of any type has the same possibility to emigrate) while all the remaining $|j|-M'$ particles reproduce according to the multivariate geometrical law defined in $\eqref{geom}$. For each $\ell\in [L]$, there is also an immigration of $A_{\ell}^{(n)}(M')+j_{\ell+1}$ new $\ell$-th type particles (here we use the convention that $j_{L+1}=0$).  
\item If $|j|< M'$, we have $Z_{n}= A^{(n)}({|j|})+\widetilde{j}.$ In this case, for each $\ell\in [L]$,  all $j_{\ell}$ particles of $\ell$-th type emigrate while $A_{\ell}^{(n)}(|j|)+j_{\ell+1}$ new particles of $\ell$-th type immigrate. 
\end{itemize}
\end{remark}
\begin{proposition}\label{ergodic}
The Markov chain $\ZZ$ is  ergodic.
\end{proposition}
\begin{proof} 
{Recall from \eqref{def.Z} that for each $n\in\Z$,  $$Z_{n}=A^{(n)}(|Z_{n-1}|)+(Z_{{n-1},2}, Z_{{n-1},3},\dots, Z_{{n-1},L},0).$$
Taking $L$ iterations,  we have that for $n\ge L$,
\begin{equation}\label{Zn} \begin{aligned}
& Z_n  = \\
& \Big(\sum_{k=1}^{L}A^{(n-k+1)}_{k}(|Z_{n-k}|), \sum_{k=2}^{L}A^{(n-k+2)}_{k}(|Z_{n-k+1}|) ,\dots,A^{(n)}_{L-1}(|Z_{n-1}|)+A^{(n-1)}_{L}(|Z_{n-2}|), A^{(n)}_{L}(|Z_{n-1}|) \Big).\end{aligned}
\end{equation}
Define $L'=\max\{\ell\in [L] : \nu(-\l)>0\}$ and $d_{\ell}=\text{Card}\big\{k\in[M]:q_k(-\ell)>0\big\}$.
Set $\mathcal{S}=\mathcal{S}_{1}\times\mathcal{S}_{2}\times\cdots\times \mathcal{S}_{L}$ where $$\mathcal{S}_{\ell}=\left\{\begin{matrix} \Z_{+}&  \text{for}& 1\le \ell \le L',\\
\{0,1,\dots, \sum_{s=\ell}^L d_{s}\} & \text{for} &  L'+1\le \l\le L.\end{matrix}\right.$$
It is evident from from \eqref{Zn} that for $n\ge L$, the support of $Z_n$ is equal to $\mathcal S$.
In particular, $\P(Z_{{L}}=j|Z_0=i)>0$ and  $\P(Z_{{L+1}}=j|Z_0=i)>0$  for any $j\in \mathcal{S}, i\in \mathbb Z_+^L$.} Hence $\ZZ$ is irreducible and aperiodic.
Using Proposition~\ref{Markovian}, we have that conditional on $\{Z_0=(0,0,\dots,0)\}$, $Z_{n-1}$ has the same distribution as $V_{0}^n$. As the process $\X$ is transient, $V_0^n$ convergences almost surely to $V_0^{\infty}=(V_{0,1}^{\infty},V_{0,2}^{\infty},\dots, V_{0,L}^{\infty})$ where $V_{0,\l}^{\infty}$ is the total number of  jumps from a site in $\N$ to site $-\l+1$. 
Hence $Z_n$ converges in law to some a.s. finite random vector $Z_{\infty}$ as $n\to\infty$.
This implies that $\ZZ$ is positive recurrent. Hence $\ZZ$ is ergodic.
\end{proof} 

\section{Critical multi-type Galton-Watson branching process with emigration}\label{sec:emigtation}
In this section, we prove a limit theorem (see Theorem \ref{lem:branching} below) for critical multi-type Galton-Watson processes with emigration. This result will be used to solve the critical case $\delta=2$ in Section \ref{sec:speed}.

\begin{definition}\label{def:emigration}
Let $N=(N_1,N_2,...,N_L)$ be a vector of $L$ deterministic positive integers and $(\psi(k,n))_{k,n\in \N}$ be a family of i.i.d. copies of a random matrices $\psi$ such that $\psi$ takes values in $\mathbb{Z}_+^{L\times L}$ and its rows are independent. Let $(U(n))_{n\ge 0}$ be a  Markov chain in $\Z_{+}^L$ defined recursively by 
$$U_j(n)=\sum_{i=1}^L \sum_{k=1}^{\varphi_i(U(n-1))} \psi_{i,j}(k,n)\quad \text{for}\ j\in [L], n\in\N$$
where $\varphi_i(s)=(s_i-N_i)\mathds{1}_{\{s_j\ge N_j, \forall j\in [L]\}}$. We call $(U(n))_{n\ge 0}$ a multi-type Galton-Watson branching process with $(N_1,N_2,\dots, N_L)$-emigration. 
\end{definition}

We can interpret the branching process $(U(n))_{n\ge0}$ defined above as a model of a population with $L$ different types of particles in which $U_i(n)$ stands for the number of particles of type $i$ in generation $n$. The number of offsprings of type $j$ produced by a particles of type $i$ has the same distribution as $\psi_{i,j}$. In generation $n$, if $U_i(n)\ge N_i$ for all $i\in [L]$ then there is an emigration of $N_i$ particles of type $i$ for $i\in [L]$, otherwise all the particles emigrate and  $U(n+1)=(0,0,\dots, 0)$. 

Let $|\cdot|$, $\|\cdot\|$ and $\langle \cdot,\cdot\rangle$ stand for the 1-norm, the Euclidean norm and the Euclidean scalar product on $\mathbb R^L$ respectively. Denote $\0=(0,0,\dots, 0)$, $\mathbf{1}=(1,1,\dots, 1)$. 

From now on, we always assume that $(U(n))_{n\ge0}$ satisfies the following condition.
\begin{assumption}\label{Asum:B} {The branching process $(U(n))_{n\ge0}$ is critical, i.e. the expected offspring matrix $\E[\psi]$ is positively regular (in the sense that there exists $n\in \N$ such that $(\E[\psi])^{n}$ is a positive matrix) and $\lambda=1$ is its largest eigenvalue in modulus. }
\end{assumption}

By Perron-Frobenius theorem, the maximal eigenvalue $1$ is simple and has the positive right and left eigenvectors $u=(u_1,\dots,u_L)$, $v=(v_1,\dots,v_L)$ which are uniquely defined such that $\langle u,v\rangle=1$ and $|u|=1$. Furthermore, other eigenvalue than $1$ is strictly smaller than 1 in modulus.

Set $$\sigma_{i,j}(k)= \E[\psi_{ki}\psi_{kj}-\delta_{i,j}\psi_{kj}],\quad \beta=\frac{1}{2}\sum_{i,j,k\in [L]}v_ku_{i}\sigma_{ij}(k)u_j\quad \text{and} \quad \theta= \frac{\langle N,u\rangle }{\beta}.$$ 

We will prove the following theorem, which is a multivariate extension of the limit theorem for critical (one-type) branching processes with emigration obtained in \cite{V1977}, \cite{V1987} and \cite{K1991}  (see also \cite{VZ1993} for a literature review). 

\begin{tcolorbox}
\begin{theorem}\label{lem:branching}
{Let $K=(K_1,K_2,\dots, K_L)$ be a deterministic vector such that $K_i\ge N_i$ for all $i\in [L]$. Assume that $U(0)=K$ a.s.} and there exists $\varepsilon>0$ such that
\begin{equation}
\label{cond:moment}\E[\psi_{i,j}^{1+{\lceil\theta\rceil}\vee (1+\varepsilon)} ]<\infty\quad \text{for all}\ \ i,j\in [L]
\end{equation}
where we denote by $\lceil\theta\rceil$ the ceiling value of $\theta$.
Then the followings hold true:
\begin{itemize}
\item[a.] There exists a constant $\varrho>0$ such that
$$\P(U(n)\neq \0)\sim \frac{\varrho}{n^{1+\theta}}\quad \text{as  } n\to\infty.$$
\item[b.] We have
$$\liminf_{n\to\infty}\frac{\E[|U_n|\ | \ U(n)\neq \0]}{n}\ge \beta.$$
\end{itemize}
\end{theorem}
\end{tcolorbox}
 For each $n\in \Z_+$ and {$s=(s_1,s_2,\dots, s_L)\in [0,1]^L$}, we set $$F(s,n):=\E\left[\prod_{\l=1}^L s_{\l}^{U_{\l}(n)}\right],\quad  f(s):=(f_1(s),\dots, f_L(s)) \text{  with  } f_i(s):=\E\left[\prod_{\l=1}^L s_{\l}^{\psi_{i,\l}}\right],$$
 which stand respectively for the multivariate probability generating functions of $U(n)$ and the random row vectors  $(\psi_{1,{\ell}})_{{\ell}=1}^L$, $(\psi_{2,{\ell}})_{{\ell}=1}^L, \dots, (\psi_{L,{\ell}})_{{\ell}=1}^L$.
 
 Let $f^k=(f^k_1,f^k_2,\dots,f^k_L)$ be the $k$-th iteration of $f$, i.e. $f^0(s)=s$ and $f^{k+1}(s)=f(f^{k}(s))$ for $k\ge0$.
 We also set $$g(s):=\prod_{\l=1}^Ls_{\l}^{-N_{\l}}\quad \text{and}\quad \gamma_n(s):=\prod_{k=1}^{n}g(f^k(s)),\quad \text{with}\quad \gamma_0(s)=1.$$

In order to prove Theorem~\ref{lem:branching}, we will need the following lemmas.   
\begin{lemma}\label{lem:gamma} Assume that there exists $\varepsilon>0$ such that
\begin{equation}\label{cond:moment2}
\E[\psi_{i,j}^{2+\varepsilon} ]<\infty,\quad\text{for all  } i,j\in [L].
\end{equation} 
Then there exists a positive constant $C$ such that 
$$\gamma_n(\0)\sim {\text{C}}{n^{\theta}} \quad \text{as}\  n\to\infty$$
and
$$\sum_{n=0}^{\infty} \gamma_{n}(\0)z^{n}\sim \frac{\Gamma(\theta+1)C}{(1-z)^{\theta+1}} \quad \text{as}\quad z\to 1^-$$
where $\Gamma$ is the Gamma function defined by $\Gamma(x)=\int_0^{\infty} u^{x-1}e^{-u}{\rm d}u$ for $x>0$.
\end{lemma}
\begin{proof}
We denote by $\mathbf 1$ the $L$-dimentional vector with all entries  equal to 1. Note that \begin{align}\label{prob.surv}\mathbf{1}-f^n(\0) \sim \frac{u}{\beta n},\quad \text{as}\ \ n\to\infty.\end{align}
(see, e.g., Corollary V.5, p.~191 in \cite{AN1972}). Set $r_n:=\left\|  \mathbf{1}-f^n(\0) -\frac{u }{\beta n}\right\|$. We first show that
\begin{align}\label{sum.r}
\sum_{n=1}^{\infty} r_n<\infty.\end{align}
Indeed, let $\mathcal Q(s)=(\mathcal Q_1(s),\mathcal Q_2(s),\dots, \mathcal Q_L(s))$ be a vector of quadratic forms with $$\mathcal Q_k(s):=\frac{1}{2}\sum_{i,j\in [L]} \sigma_{ij}(k)s_i s_j$$ and set
$$a(s):=\left\langle v, \mathcal Q\left(\frac{\mathbf 1-s}{\langle v,\mathbf 1-s\rangle}\right) \right\rangle,\quad d(s):=\frac{1}{\langle v,\mathbf 1-s\rangle}+a(s)-\frac{1}{\langle v,\mathbf 1-f(s)\rangle}.$$
We have
\begin{align}\label{diff}{
\frac{1}{\langle v,\mathbf 1-f^n(s)\rangle}-\frac{1}{\langle v,\mathbf 1-s\rangle}=\sum_{k=0}^{n-1}a(f^k(s))-\sum_{k=0}^{n-1}d(f^k(s)).}\end{align}
In virtue of Taylor's expansion, we have $\mathbf 1-f(s)=(\E[\psi]-H(s))\cdot (\mathbf 1-s)$
with $H(s)=O(\|\mathbf 1-s\|)$. It follows that
\begin{align}\label{prod}
\mathbf 1-{f^n(\0)}=\left(\E[\psi]-H(f^{n-1}(\0))\right)\dots \left(\E[\psi]-H(f(\0))\right)\left(\E[\psi]-H(\0)\right)\cdot\mathbf 1.\end{align}
{Using \eqref{prob.surv},} we note that $\|H(f^n(\0))\|=O(\|\mathbf 1-f^n(\0)\|)=O(1/n)$. Moreover, $\E[\psi]^n=uv^T+O(|\lambda|^{n})$ where $\lambda$ (with $|\lambda|<1$) is the eigenvalue of $\E[\psi]$ with the second largest modulus. In what follows, we denote by $\textsf{Cst}$ a positive constant but its value may vary from line to line. Applying inequality (4.11) in \cite{JS1967} to \eqref{prod}, we deduce that
\begin{equation}\label{eq.fn}\left\|\frac{\mathbf 1-f^n(\0)}{\langle v,\mathbf  1-f^n(\0)\rangle}-u \right\|\le \frac{\textsf{Cst}}{n}.
\end{equation}
Since $\mathcal Q(s)$ is Lipschitz, we thus have $$\left\|a(f^n(\0))-\langle v,\mathcal Q(u)\rangle\right\| \le \|v\|.\left\|\mathcal Q\left(\frac{\mathbf  1-f^n(\0)}{\langle v,\mathbf  1-f^n(\0)\rangle}\right)-\mathcal Q(u) \right\|\le \frac{\textsf{Cst}}{n}.$$
As a result, we have
\begin{align}\label{eq.a}{\sum_{k=0}^{n-1}a(f^k(\0))=\langle v,\mathcal Q(u)\rangle n +O(\log(n))=\beta n+O(\log(n)).}
\end{align}
W.l.o.g., we assume that $\varepsilon\in(0,1)$ (which satisfies \eqref{cond:moment2}). By Taylor's expansion, there exists a vector function $\mathcal E(t,s)$ such that
$\mathbf 1-f(s)=\E[\psi](\mathbf 1-s)-\mathcal Q(\mathbf 1-s)+\mathcal E(\mathbf 1-s,\mathbf 1-s)$
where we note that \begin{equation}\label{vfunc.E}\mathcal \|\mathcal E(t,s)\|\le \textsf{Cst}  \|t\|^{\varepsilon}\| s\|^2.\end{equation} 
According to inequality (4.42) in \cite{JS1967}, we have  
\begin{align}\label{ineq.d1}
 -\langle v, \mathbf 1- f^n(\0)\rangle \left\langle v, \mathcal Q\left( \frac{\mathbf 1-f^n(\0)}{\langle v, \mathbf 1-f^n(\0) \rangle} \right)\right \rangle^2 \le
d(f^n(\0))\le \left\langle v, \mathcal E\left(\mathbf 1-f^n(\0), \frac{\mathbf 1-f^n(\0)}{\langle v, \mathbf 1-f^n(\0) \rangle} \right)\right \rangle.
\end{align}
By reason of \eqref{eq.fn}, we notice that
\begin{equation}\label{lim.fn}
\frac{\mathbf 1-f^n(\0)}{\langle v, \mathbf 1-f^n(\0) \rangle}\to u\quad \text{as} \ \ n\to\infty.
\end{equation}
Combining \eqref{ineq.d1} with \eqref{vfunc.E}-\eqref{lim.fn} and using Cauchy-Schwarz inequality, we obtain
\begin{align}\label{ineq.d2}
d(f^n(\0))\ge  - \textsf{Cst} \| \mathbf 1- f^n(\0)\|= O(n^{-1}) \quad\text{and}\quad d(f^n(\0))\le \textsf{Cst}\|\mathbf 1-f^n(\0)\|^{\varepsilon}=O(n^{-\varepsilon}).
\end{align}
Combining \eqref{diff} with \eqref{eq.a} and \eqref{ineq.d2}, we get
$${\langle v, \mathbf 1- f^n(\0)\rangle =\frac{1}{|v|^{-1}+ \sum_{k=0}^{n-1}a(f^k(\0))-\sum_{k=0}^{n-1}d(f^k(\0)) } =\frac{1}{\beta n } + O(n^{-1-\varepsilon}).} $$
Consequently,
\begin{align}\label{fn0} \begin{aligned}\mathbf 1- f^n(\0) =\langle v, \mathbf 1- f^n(\0) \rangle.\frac{\mathbf 1- f^n(\0)}{ \langle v, \mathbf 1- f^n(\0) \rangle}  & = \left( \frac{1}{\beta n } + O(n^{-1-\varepsilon})\right)(u+O(n^{-1}))\\ & = \frac{u}{\beta n}+O(n^{-1-\varepsilon}).\end{aligned}\end{align}
Hence $r_n=O(n^{-1-\varepsilon})$ and \eqref{sum.r} is thus proved. On the other hand, by Taylor's expansion, we have 
\begin{align}\label{gf0}g(f^k(\0)) 
& = 1 + \langle N , \mathbf{1}-f^k(\0) \rangle  + O(\|\mathbf{1}-f^k(\0)\|^2) .\end{align}
Thus
\begin{align*}\gamma_n(\0) & =  \prod_{k=1}^{n}\left( 1+ \langle N,\mathbf{1}-f^k(\0) \rangle + O(|\mathbf{1}-f^k(\0)|^2) \right) \\
&
 \sim \textsf{Cst}\cdot \exp\left(\sum_{k=1}^{n} \left[\langle N, \mathbf{1}-f^k(\0) \rangle + O(k^{-2}) \right]\right)\\
& \sim \textsf{Cst}\cdot \exp\left(\frac{\langle N,u \rangle}{\beta}\sum_{k=1}^{n}\left(\frac{1}{k}+O(k^{-(1+\varepsilon)})\right)\right). 
\end{align*}
Since $\sum_{k=1}^n1/k=  \log(n)+O(1)$ as $n\to\infty$ and  $\theta=\langle N,u \rangle/{\beta}$, we obtain that $\gamma_n(\0)\sim C n^{\theta}$ for some positive constant $C$. Furthermore, by Hardy–Littlewood tauberian theorem for power series (see e.g. Theorem 5, Section XIII.5, p. 447 in \cite{F1971}), we deduce that
$$\sum_{n=0}^{\infty} \gamma_{n}(\0)z^{n}\sim \frac{\Gamma(\theta+1)C}{(1-z)^{\theta+1}} \quad \text{as}\quad z\to 1^-.$$
\end{proof}
In what follows, for $x,y \in \Z_+^{L}$ we write $x\succeq y$  if $x_{i}\ge y_i$ for all $i\in [L]$, otherwise we write $x\nsucceq y$. Set $\mathcal{S}(N)=\{r\in {\Z_{+}^L\setminus\{\mathbf{0}\}}: r\nsucceq  N \}$.
For each $r\in \mathbb Z_{+}^L$ and $s=(s_1,s_2,\dots, s_L)\in \mathbb \R^L$ such that $s_{\ell}\neq0$ for all $\l\in[L]$, define
$$H_r(s):=\left(\prod_{\l=1}^Ls_{\l}^{r_{\l}-N_{\l}}-1\right)\mathds{1}_{\{ r \in \mathcal{S}(N)\}}.$$
For each $n\in \Z_+$ and $z\in [0,1)$, set
$$\mu_{n}:=\P(U(n)\neq \0)=1-F(\0,n)\quad\text{and}\quad Q(z):=\sum_{n=0}^{\infty} \mu_{n} z^n. $$
\begin{lemma}\label{lem:Q} The generating function of $(\mu_n)_{n\ge0}$ is given by
$$Q(z) =\frac{B(z)}{D(z)},$$
in which we define \begin{align*} B(z)& :=\sum_{n=0}^{\infty}\left(1-F(f^n(\0),0)+\sum_{k=1}^{\infty}\E\left[H_{U(k-1)}(f^{n+1}(\0))  \right]z^k\right){\gamma_{n}(\mathbf{0})}z^n\quad \text{and}\\
D(z)&:=(1-z)\sum_{n=0}^{\infty} \gamma_{n}(\0)z^{n}.
\end{align*}
\end{lemma}
 \begin{proof}
From the definition of $U(n)$, we have
\begin{align*}F(s,n)& =\E\Big[\prod_{i=1}^L \prod_{k=1}^{\varphi_i(U_i(n-1))}\prod_{\l=1}^L s_{\l}^{\psi_{i,\l}(k,n)} \Big]= \E\Big[\prod_{i=1}^L (f_i(s))^{\varphi_i(U_i(n-1))}  \Big]\\
& = \sum_{r\in \Z_+^L, r \succeq  N} \P( U({n-1})=r)\prod_{i=1}^L (f_i(s))^{r_i-N_i} +  \P( U({n-1})=\0)+\sum_{r\in \mathcal{S}(N)} \P( U({n-1})=r)
\\&  = F(f(s),n-1)g(f(s))-F(\0,n-1)(g(f(s))-1) -\sum_{r\in \Z_+^L}\P(U(n-1)=r) H_r(f(s)).
\end{align*}
Consequently,
\begin{align}\label{eqn.F} F(s,n)& =F(f^n(s),0)\gamma_{n}(s)-\sum_{k=1}^n F(\0, n-k)(\gamma_k(s)-\gamma_{k-1}(s))
\\
& \nonumber-\sum_{r\in \Z_+^L  }\sum_{k=1}^{n}\P(U(n-k)=r)H_{r}(f^k(s))\gamma_{k-1}(s).\end{align}
Note that {$\mu_n=1-F(\0,n)$}. Substituting $s=\0$ into \eqref{eqn.F}, we obtain
\begin{align*}\mu_n +\sum_{k=1}^n \mu_{n-k}(\gamma_{k}(\0)-\gamma_{k-1}(\0))&=\left(1-F(f^n(\0),0)\right)\gamma_{n}(\0)\\
&+\sum_{r\in \Z_+^L  }\sum_{k=1}^n \P(U(n-k)=r)H_{r}(f^k(\0))\gamma_{k-1}(\0).\end{align*}
Multiplying both sides of the above equation by $z^n$ and summing over all $n\ge0$, we get 
\begin{align*}{(1-z)\Big(\sum_{n=0}^{\infty}\gamma_{n}(\0)z^{n}\Big)\Big( \sum_{n=1}^{\infty}\mu_nz^{n}\Big) }& = \sum_{n=0}^{\infty}\big(1-F(f^n(\0),0)\big)\gamma_{n}(\0)z^{n} 
\\ 
& +\sum_{r\in \Z_+^L  }\sum_{k=1}^{\infty} \P(U(k-1)=r)z^{k}\sum_{n=0}^{\infty}H_{r}(f^{n+1}(\0))\gamma_{n}(\0)z^{n}.
\end{align*}
This ends the proof of the lemma.
\end{proof} 

Let $(p_{n})_{n\ge0}$  be a non-negative sequence such that its generating function $P(z)=\sum_{n=0}^{\infty}p_nz^n$ has radius of converge 1. For each $k\in \Z_+$, we define the sequence $(p_n^{(k)})_{n\ge 0}$ recursively by $$p_n^{(0)}= p_n\quad \text{and}\quad { p_n^{(k+1)}=\sum_{j=n+1}^{\infty}p^{(k)}_j \quad \text{for} \ k\in\Z_+}.$$ 
Denote by $P^{(k)}(z)=\sum_{n=0}^{\infty}p_{n}^{(k)}z^n$ the generating function of $(p_n^{(k)})_{n\ge0}$. Assume that
$P^{(j)}(1^-)<\infty$ for all $0\le j\le k$. By Abel's theorem, we notice that
$P^{(j)}(1^-)=P^{(j)}(1)=\sum_{n=0}^{\infty}p_{n}^{(j)}$ and 
$$P^{(j+1)}(z)=\frac{P^{(j)}(1^-)-P^{(j)}(z)}{1-z}\quad \text{for }z\in (-1,1), 0\le j\le k.$$

{The next lemma is derived directly from Corollary 2 and Lemma 5 in \cite{V1977}.
\begin{lemma}\label{iter.sum} Let $(p_n)_{n\ge0}$ be a non-negative non-increasing sequence and $p_n\to 0$ as $n\to\infty$.\begin{itemize}
\item[i.] Assume that $p^{(1)}_n\sim n^{-\alpha}$ as $n\to\infty$ for some $\alpha>0$. Then there exists a constant $\varpi>0$ such that $p_n\sim \varpi n^{-\alpha-1}$ as $n\to\infty$.
\item[ii.] Assume that $\sum_{k=1}^{n}kp_k\sim 
n^{2-\alpha}$ as $n\to\infty$ for some $0<\alpha<2$. Then $p_n\sim (2-\alpha)n^{-\alpha}$ as $n\to\infty$.
\end{itemize}
\end{lemma}}

For $r\in\Z_{+}^L$, $n\in\N$ and $k\in\Z_+$ set $\pi_{r,n}:=\P(U(n-1)=r),$   $$\Pi_r(z):=\sum_{n=1}^{\infty} \pi_{r,n} z^n\quad\text{and}\quad  \Pi_r^{(k)}(z):=\sum_{n=1}^{\infty} \pi_{r,n}^{(k)} z^n.$$
Recall that $\mathcal{S}(N)=\{r=(r_1,r_2,\dots,r_L)\in\Z_+^L\setminus\{\0\}: \exists \ell\in[L], r_{\ell}<N_{\ell} \}$.
\begin{lemma}\label{lem.Pi}
We have $\sum_{r\in \mathcal{S}(N)}\Pi_r(1^{-}) \le 1.$
Furthermore, if $Q^{(k)}(1^{-})<\infty$ for some $k\ge0$ then $\sum_{r\in \mathcal{S}(N)}\Pi_r^{(k+1)}(1^{-})<\infty.$
\end{lemma}

\begin{proof}
Define $\tau=\inf\{n\ge 0: U(n+1)=\0 \}$. For each $r\in \mathcal{S}(N)$, we have $$\{U(n)=r\}=\{U(n)=r, U(n+1)=\0\}=\{U(\tau)=r, \tau=n\}$$
yielding that \begin{align}\label{eq.Pi}\sum_{r\in \mathcal{S}(N)}\Pi_r(1^{-})=\sum_{r\in \mathcal{S}(N)}\sum_{n=0}^{\infty} \P(U(n)=r)=\sum_{r\in \mathcal{S}(N)}\P(U(\tau)=r)\le 1.\end{align}
We also have
\begin{align}\label{inq.pi1}
\sum_{r\in \mathcal{S}(N)} \pi_{r,n}^{(1)}&=\sum_{r\in \mathcal{S}(N)}\sum_{m=n+1}^{\infty} \pi_{r,m} = \sum_{r\in \mathcal{S}(N)}\sum_{m=n+1}^{\infty}  \P(U(m-1)=r)\\
\nonumber&=\sum_{r\in \mathcal{S}(N)}\P\left(U(\tau)=r, \tau\ge n\right)\le \P(\tau\ge n)=\mu_{n}.
\end{align}
By induction, we obtain that $\sum_{r\in \mathcal{S}(N)}\pi_{r,n}^{(k+1)}\le \mu_{n}^{(k)}$ for all $k\ge0$. Hence
$$\sum_{r\in \mathcal{S}(N)}\Pi_r^{(k+1)}(1^{-})=\sum_{r\in \mathcal{S}(N)}\sum_{n=0}^{\infty} \pi_{r,n}^{(k+1)}\le \sum_{n=0}^{\infty} \mu_{n}^{(k)}=Q^{(k)}(1^{-}).$$
This ends the proof of the lemma.
\end{proof}

{In what follows, we set $\widetilde{\theta}:=\lceil \theta \rceil$, which is the ceiling value of $\theta$. We denote by $\textsf{Cst}$ a generic non-zero constant and its value may vary from line to line.}

\begin{lemma}\label{lemmaDB} Assume that the condition \eqref{cond:moment} is fulfilled. Then:

i. For each $ k\in\{0,1,\dots, \widetilde{\theta}-1\}$, there exist power series $R_D^{[k]}(z), D_m^{[k]}(z)$ and non-zero constants $d^{[k]}_m$, $m\in \{k+1,\dots, \widetilde{\theta}\}$  such that as $z\to 1^{-}$,
\begin{align}
\label{RD.gen} 
& {R_D^{[k]}(z)\sim \left\{\begin{matrix}{\sf Cst}\cdot \ln\big(\frac{1}{1-z}\big) & \text{if} \ \theta\in{\N},\\
{\sf Cst}& \text{otherwise},
\end{matrix}\right.}\\
\label{Dm.gen} 
&D_m^{[k]}(z)\sim \Gamma(\theta-m+1)d^{[k]}_m.(1-z)^{-(\theta-m+1)}\quad \text{and}\\
\label{D.gen}
&(1-z)^kD(z)=\sum_{m=k+1}^{\widetilde{\theta}}D_m^{[k]}(z)+R_D^{[k]}(z).
\end{align}
ii. There exist power series $R_B(z), B_m(z)$ and $b_m(z)$ with $m\in\{1,2,\dots,\widetilde{\theta}\}$ such that $|b_m(1^-)|<\infty$,
\begin{align}
\label{RB.gen}
&{R_B(z)\sim 
\left\{\begin{matrix} {\sf Cst}\cdot\ln\big(\frac{1}{1-z}\big), & \text{if } \theta\in\N,\\ {\sf Cst} &\text{otherwise,}\end{matrix}\right.}\\
\label{Bm.gen}
& B_m(z)\sim \Gamma(\theta-m+1)b_m(z). (1-z)^{-(\theta-m+1)}\quad \text{and}\\
\label{B.gen}
&B(z)=\sum_{m=1}^{\widetilde{\theta}}B_{m}(z)+R_B(z)
\end{align}
as $z\to 1^-$.
\end{lemma}
\begin{proof}
For $k\in \N$, $r\in \Z^L\setminus\{\0\}$ and $j=(j_1,\dots, j_L)\in \Z_+^L$ with $|j|\le \widetilde{\theta}+1$, we set 
\begin{equation}\label{coef.c}
c_{j,r,k}:=\frac{1}{j_{1}!\dots j_{L}!}\,\left.\frac{\partial^{j_{1}+\cdots +j_{d}}\prod_{\l=1}^L(f_{\l}^{k}(s))^{r_{\l}}}{\partial s_{1}^{j_{1}}\cdots \partial s_{d}^{j_{d}}}\right|_{s=\mathbf 1}\end{equation}
which is well-defined thanks to the condition \eqref{cond:moment}. Using Taylor's expansion, we have that for $r=(r_1,\dots, r_L)\in\Z^L\setminus\{\0\}$, $n\ge \widetilde{\theta}$ and $k\ge 1$, 
\begin{align}\label{taylor} \prod_{\l=1}^L(f_{\l}^{n-\widetilde{\theta}+k}(\0))^{r_{\l}} & =1 + \sum _{1\le |j|\le \widetilde{\theta}}(-1)^{|j|} c_{j,r, k }  \prod_{\l=1}^L (1-f_{\l}^{n-\widetilde{\theta}}(\0))^{j_{\l}}+{\varrho_{n,r,k}}\\
{\nonumber\text{with}  \quad \varrho_{n,r,k} } & { \sim (-1)^{\widetilde{\theta}+1}\sum_{|j|=\widetilde{\theta}+1}c_{j,r,k}\prod_{\l=1}^L (1-f_{\l}^{n-\widetilde{\theta}}(\0))^{j_{\l}}\text{  as  } n\to\infty.}\end{align}
{By reason of \eqref{prob.surv}, we notice that
\begin{align}\label{Res}\varrho_{n,r,k}\sim (-1/\beta)^{\widetilde{\theta}+1} n^{-\widetilde{\theta}-1}\sum_{|j|=\widetilde{\theta}+1}c_{j,r,k}\prod_{\l=1}^L u_\ell^{j_\ell}.
\end{align}
Recall from Lemma \ref{lem:gamma} that $\gamma_n(\0)\sim C n^{\theta}.$ Using Hardy–Littlewood tauberian theorem for power series, we thus have that as $z\to 1^{-}$,
\begin{align}\label{sum.R}\sum_{n=\widetilde{\theta}}^{\infty} \varrho_{n,r,k} \gamma_{n}(\0) z^{n} \sim \left\{\begin{matrix} \textsf{Cst} \cdot\ln\big(\frac{1}{1-z}\big), & \text{if } \theta\in\N,\\ \textsf{Cst} &\text{otherwise}.
\end{matrix} \right.\end{align}}

{\underline{Part i.} Notice that for $n\in\N$, $\gamma_{n}(\0)-\gamma_{n-1}(\0)=\gamma_{n}(\0)\big(1-\prod_{\ell=1}^L(f^n_{\l}(\0))^{N_{\l}}\big)$ and thus
\begin{align}\label{D.expand}
\nonumber D(z)&=(1-z)\sum_{n=0}^{\infty}\gamma_n(\0)z^n=\gamma_0(\0)+\sum_{n=1}^{\infty}(\gamma_n(\0)-\gamma_{n-1}(\0))z^n\\
&=\sum_{n=0}^{\infty}\left(1-\prod_{\l=1}^L(f^{n}_{\l}(\0))^{N_{\l}} \right)\gamma_n(\0)z^n.
\end{align}
By induction, we easily obtain that for $ k\in\{0,1,\dots, \widetilde{\theta}-1\},$
\begin{align}\label{Dk}
   (1-z)^kD(z) = (1-z)^{k+1}\sum_{n=0}^{\infty}\gamma_n(\0)z^n = \sum_{n=0}^{\infty}\prod_{m=0}^{k}\left(1-\prod_{\l=1}^L(f^{n-m}_{\l}(\0))^{N_{\l}} \right)\gamma_n(\0)z^n
    \end{align}
where we use the convention that $f^{h}(\0)=\0$ for $h\le 0$.}
Using \eqref{taylor} and \eqref{Dk}, we have that for $k\in\{0, 1,\dots, \widetilde{\theta}-1\},$
\begin{align*}(1-z)^kD(z) 
& = \sum_{n=\widetilde{\theta}}^{\infty}\left[\sum _{k+1\le|j|\le \widetilde{\theta}+1 } \chi_{j}^{[k]} \prod_{\l=1}^L (1-f_{\l}^{n- \widetilde{\theta} }(\0))^{j_{\l}} + o(\|\mathbf 1- f^{n- \widetilde{\theta}}(\0) \|^{\widetilde{\theta}+1})\right]\gamma_{n}(\0)z^{n}\\
& + \sum_{n=0}^{\widetilde{\theta}-1}\prod_{m=0}^{k}\left(1-\prod_{\l=1}^L(f^{n-m}_{\l}(\0))^{N_{\l}} \right)\gamma_n(\0)z^n
\end{align*}
where for $k\in \Z_{+}$ and $j\in \Z_{+}^L$, we set
$$\chi_j^{[k]}:=(-1)^{|j|+k+1}\sum_{\substack{ j^{(m)}\in \Z_{+}^L\setminus\{\0\},\\j^{(0)}+\dots+j^{(k)}=j }}\prod_{m=0}^{k}c_{j^{(m)},N,\widetilde{\theta}-m }.$$ 
In virtue of Lemma~\ref{lem:gamma} and \eqref{prob.surv}, we note that \begin{equation}\label{asym.term}
\gamma_{n}(\0)\prod_{\l=1}^L (1-f_{\l}^{n- \widetilde{\theta} }(\0))^{j_{\l}}\sim C\beta^{-|j|}\left(\prod_{\l=1}^L u_{\l}^{j_{\l}}\right) n^{\theta-|j|} \quad\text{as}\ n\to\infty.\end{equation}
For $1\le m\le \widetilde{\theta}+1$ and $0\le k\le \widetilde{\theta}-1$, set
$$d_m^{[k]}:= C\beta^{-m}\sum _{|j|=m } \chi_j^{[k]}\prod_{\l=1}^L u_{\l}^{j_{\l}} \quad\text{and}\quad D_{m,n}^{[k]}:=\gamma_{n}(\0)\sum _{|j|=m } \chi_j^{[k]} \prod_{\l=1}^L (1-f_{\l}^{n- \widetilde{\theta} }(\0))^{j_{\l}}.$$ 
Note that $D_{m,n}^{[k]}\sim d_m^{[k]} n^{\theta-m}$ as $n\to\infty$. For $m\in\{k+1,\dots, \widetilde{\theta}\}$, define
$D_m^{[k]}(z):=\sum_{n=\widetilde{\theta}}^{\infty} D_{m,n}^{[k]}z^{n}$. 
In view of Hardy–Littlewood tauberian theorem for power series, we thus obtain \eqref{RD.gen}, \eqref{Dm.gen} and \eqref{D.gen}.

\underline{Part ii.} Recall that $$B(z)=\sum_{n=0}^{\infty}\left[1-F(f^n(\0),0)+\sum_{r\in \mathcal{S}(N)} H_{r}(f^{n+1}(\0)) \Pi_r(z)\right]\gamma_{n}(\0)z^{n}$$
in which we have $$F(f^n(\0),0)=\prod_{\l=1}^L(f^n_{\l}(\0))^{K_{\l}} \quad\text{and}\quad H_{r}(f^{n+1}(\0)) =\prod_{\l=1}^L(f^{n+1}_{\l}(\0))^{r_{\l}-N_{\l}}-1 \ \text{ for } r\in\mathcal{S}(N).$$
Using \eqref{taylor}, we obtain that
\begin{align*}B(z) & =  \sum_{n=\widetilde{\theta}}^{\infty}\left[\sum _{0<|j|\le \widetilde{\theta} } a_{j}(z) \prod_{\l=1}^L (1-f_{\l}^{n- \widetilde{\theta} }(\0))^{j_{\l}}+{ \varrho_n^B(z)}\right]\gamma_{n}(\0)z^{n}\\
&+ \sum_{n=0}^{\widetilde{\theta}-1}\left[1-F(f^n(\0),0)  +\sum_{r\in \mathcal{S}(N)}H_{r}(f^{n+1}(\0)) \Pi_r(z)\right]\gamma_{n}(\0)z^{n}
\end{align*}
where we set 
\begin{align*}a_j(z)&:= (-1)^{|j|-1}c_{j,K,\widetilde{\theta}}+(-1)^{|j|} \sum_{r\in \mathcal{S}(N)}  {c_{j,r-N, \widetilde{\theta}+1}}\Pi_r(z),\\
\varrho^B_n(z)&:= - \varrho_{n,K,\widetilde{\theta}}+\sum_{r\in \mathcal{S}(N)} \varrho_{n,r-N, \widetilde{\theta}+1}\Pi_r(z).\end{align*}
Here we notice that for $s\in\R^L_+$ with $|s|< 1$ we have
$$\Big| \sum_{r\in \mathcal{S}(N)} H_{r}\big(f_{\l}^{\widetilde{\theta}+1}(s)\big) \Pi_r(1^-)\Big| \le \left(\prod_{\ell=1}^L\big( f_{\l}^{\widetilde{\theta}+1}(\0)\big)^{-N_{\ell}}+1\right)\sum_{r\in \mathcal{S}(N)} \Pi_r(1^-)<\infty.$$
It follows that
\begin{equation}\label{sumable}\Big|\sum_{r\in \mathcal{S}(N)}  {c_{j,r-N, \widetilde{\theta}+1}}\Pi_r(1^-)\Big|<\infty\quad\text{and thus}\quad |a_j(1^-)|<\infty \quad \text{for all }1\le |j|\le \widetilde{\theta}+1.
\end{equation}
By \eqref{Res} and \eqref{sumable}, we also note that $|\varrho_n^B(1^-)|<\infty$. Furthermore, in virtue of \eqref{sum.R}, we notice that as $z\to 1^{-}$,
\begin{align}\label{Res.B}\sum_{n=\widetilde{\theta}}^{\infty}\varrho_n^{B}(z)\gamma_{n}(\0)z^{n}=\left\{\begin{matrix} \textsf{Cst}\cdot\ln\big(\frac{1}{1-z}\big), & \text{if } \theta\in\N,\\ \textsf{Cst} &\text{otherwise.} 
\end{matrix} \right.\end{align}
For $m\in\{1,2,\dots, \widetilde{\theta}\}$, set
\begin{align*} b_m(z)&:= (-1)^{m}C\beta^{-m} \sum_{|j|=m}  a_j(z)\prod_{\l=1}^L u_{\l}^{j_{\l}} \quad \text{and}\\ B_{m,n}(z)&:=\gamma_{n}(\0)\sum_{|j|=m} a_j(z) \prod_{\l=1}^L (1-f_{\l}^{n- \widetilde{\theta} }(\0))^{j_{\l}}\sim b_m(z)n^{\theta-m}\quad\text{as} \ n\to\infty\end{align*}
and define $B_{m}(z):=\sum_{n=\widetilde{\theta}}^{\infty}B_{m,n}(z)z^{n}.$
By \eqref{sumable}, we note that $|b_m(1^-)|<\infty$. In view of Hardy–Littlewood tauberian theorem and \eqref{Res.B}, we obtain \eqref{RB.gen}, \eqref{Bm.gen} and \eqref{B.gen}. 
\end{proof}

We notice that
$$\sum_{n=0}^{\infty}\mu_{n}=Q(1^{-}) =\lim_{z\to 1^{-}}\frac{B(z)(1-z)^{\theta}}{D(z)(1-z)^{\theta}}=\frac{b_1(1^-)}{d_1^{[0]}}<\infty.$$

\begin{proof}[Proof of Theorem~\ref{lem:branching}]
{\underline{Part a.}} We adopt an idea by Vatutin (see \cite{V1977}) as follows. Recall that $\widetilde{\theta}:=\lceil \theta \rceil$ is the ceiling value of $\theta$. We will prove that for all $k\in\{0,1,\dots, \widetilde{\theta}-1\}$, 
\begin{align}\label{sum.q}
Q^{(k)}(1^-)=\sum_{n=0}^{\infty} \mu^{(k)}_n <\infty
\end{align}
and as $z\to 1^-$,
\begin{align}\label{sum.q2}
&Q^{(\widetilde{\theta})}(z)=\sum_{n=0}^{\infty} \mu_{n}^{(\widetilde{\theta})}z^n  \sim{\left\{\begin{matrix}{\sf Cst}\cdot \ln\big(\frac{1}{1-z}\big) & \text{if } \theta\in\N,\\
 {\sf Cst}\cdot (1-z)^{-(\widetilde{\theta}-\theta)} & \text{otherwise,}\end{matrix}\right.}\\
&\label{dQtheta}{\frac{\text{d} }{\text{d} z } Q^{(\widetilde{\theta})}(z)\sim \textsf{Cst}\cdot (1-z)^{-1}\quad \text{if } \theta\in \N.}
\end{align}

{We first consider the case when $\theta$ is not an positive integer. By Hardy–Littlewood tauberian theorem for power series, it follows from \eqref{sum.q2} that $\mu_{n}^{(\widetilde{\theta})}\sim {\sf Cst}\cdot n^{\widetilde{\theta}-\theta-1}$ as $n\to\infty$.  If $\theta\in \mathbb{N}$ then it follows from \eqref{dQtheta} and Hardy–Littlewood tauberian theorem for power series that 
$\sum_{k=1}^n k \mu^{(\theta)}_k\sim {\sf Cst} \cdot n$ as $n\to\infty$. In the latter case, by Lemma \ref{iter.sum}.ii, we obtain that $\mu^{(\theta)}_n\sim \textsf{Cst}\cdot n^{-1}$ as $n\to\infty$. In both cases, using Lemma \ref{iter.sum}.i, we deduce that $\mu_{n}\sim \varrho n^{-\theta-1}$ for some constant $\varrho>0$ as $n\to\infty$.
}

Hence, to finish {the proof of Theorem~\ref{lem:branching}(a)}, we only have to verify \eqref{sum.q}, \eqref{sum.q2} and \eqref{dQtheta}.

Set $B^{[0]}(z)=B(z)$ and $B^{[k]}(z)=Q^{(k-1)}(1^-)(1-z)^{k-1}D(z)-B^{[k-1]}(z)$ for $1\le k\le \widetilde{\theta}.$ {By induction on $k$}, we notice that if $Q^{(k-1)}(1^-)<\infty$ then 
$$Q^{(k)}(z)=\sum_{n=0}^{\infty}\mu_{n}^{(k)}z^n  =\frac{Q^{(k-1)}(1^-)-Q^{(k-1)}(z)}{1-z}=\frac{B^{[k]}(z)}{(1-z)^kD(z)}.$$

Assume that up to some $k\in\{1,2,\dots \widetilde{\theta}\}$, the power series $R^{[k-1]}_B(z)$, $B_m^{[k-1]}(z)$, and $b_m^{[k-1]}(z)$ are defined  for all $m\in \{k,k+1,\dots, \widetilde{\theta}\}$ such that $|b_m^{[k-1]}(1^-)|<\infty$ and as $z\to 1^{-}$,
\begin{align*}
{R_{B}^{[k-1]}(z)
}& \sim
\left\{\begin{matrix} \textsf{Cst}\cdot\ln\big(\frac{1}{1-z}\big) & \text{if } \theta\in\N,\\ \textsf{Cst} &\text{otherwise,}\end{matrix}\right.\\
B_m^{[k-1]}(z)& \sim \Gamma(\theta-m+1)b^{[k-1]}_m(z). (1-z)^{-(\theta-m+1)},  \\
B^{[k-1]}(z) & =\sum_{m=k}^{\widetilde{\theta}}B_m^{[k-1]}(z) +R_{B}^{[k-1]}(z)
\end{align*}
yielding that $Q^{(k-1)}(1^-)={b^{[k-1]}_k(1^-)}/{d^{[k-1]}_k}$ is finite. The above statement holds for $k=1$ thanks to Lemma \ref{lemmaDB}. We next prove that this statement also holds true when replacing $k$ by $k+1$. Indeed, by reason of Lemma \ref{lemmaDB}, notice that as $z\to1^{-}$,
\begin{align*}B^{[k]}(z)=
 \sum_{m=k}^{\widetilde{\theta}}\left(Q^{(k-1)}(1^-)D_m^{[k-1]}(z)- B_m^{[k-1]}(z)\right)+ Q^{(k-1)}(1^-)R^{[k-1]}_D(z)- R_B^{[k-1]}(z).\end{align*}
We also have that as $z\to 1^{-}$,
\begin{align*}
R_B^{[k]}(z):= Q^{(k-1)}(1^-)R^{[k-1]}_D(z)- R_B^{[k-1]}(z)& \sim \left\{\begin{matrix} \textsf{Cst}\cdot\ln\big(\frac{1}{1-z}\big) & \text{if } \theta\in\N,\\ \textsf{Cst} &\text{otherwise,} \end{matrix}\right. \\
 \text{and}\quad Q^{(k-1)}(1^-)D_k^{[k-1]}(z)- B_k^{[k-1]}(z) & \sim\Gamma(\theta-k+1) \left(b_{k}^{[k-1]}(1^-)-b_{k}^{[k-1]}(z)\right)(1-z)^{\theta}\\
&= \Gamma(\theta-k+1) \widehat{b}_k(z)(1-z)^{-(\theta-1)},
\end{align*}
where we set $\widehat{b}_k(z):=(b_{k}^{[k-1]}(1^-)-b_{k}^{[k-1]}(z))/(1-z)$.
For $m\in\{k+1, \dots,\widetilde{\theta} \}$, set
$$B^{[k]}_{m}(z):=\left\{\begin{array}{ll}\displaystyle
\sum_{n\in\{k,k+1\}}\left(Q^{(k-1)}(1^-)D_n^{[k-1]}(z)- B_n^{[k-1]}(z)\right) & \text{if}\ \ m= k+1,\\
\displaystyle Q^{(k-1)}(1^-)D_{m}^{[k-1]}(z)-B^{[k-1]}_{m}(z) & \text{if}\ \ k+2\le m\le \widetilde{\theta}\end{array}\right.$$
and 
$$b^{[k]}_{m}(z):=\left\{\begin{array}{ll}\displaystyle (\theta-k)\widehat{b}_k(z) + Q^{(k-1)}(1^-)d_{k+1}^{[k-1]}-b_{k+1}^{[k-1]}(z)
  & \text{if}\ m= k+1,\\
\displaystyle Q^{(k-1)}(1^-)d_{m}^{[k-1]}-b_{m}^{[k-1]}(z)  & \text{if}\ k+2\le m\le \widetilde{\theta}.\end{array}\right.$$
By the recurrence relation of {$b^{[k]}_m(z)$} and Lemma~\ref{lem.Pi}, we note that {$|b^{[k]}_m(1^-)|<\infty$}.  Therefore, as $z\to 1^-$,
\begin{align*}
 & B^{[k]}(z) =\sum_{m=k+1}^{\widetilde{\theta}}B_m^{[k]}(z) +R_B^{[k]}(z)\quad \text{with}\quad R_B^{[k]}(z) \sim \left\{\begin{matrix} \textsf{Cst}\cdot\ln\big(\frac{1}{1-z}\big) & \text{if } \theta\in\N,\\ \textsf{Cst} &\text{otherwise} \end{matrix}\right.\\
\text{and}\quad
& B_m^{[k]}(z) \sim \Gamma(\theta-m+1)b^{[k]}_m(z). (1-z)^{-(\theta-m+1)}  
\end{align*}
and thus $Q^{(k)}(1^-)={b_{k+1}^{[k]}(1^-)}/{d^{[k]}_{k+1}}$ is finite. By the principle of mathematical induction, we deduce that \eqref{sum.q} holds true for all $k\in \{0,1,\dots, \widetilde{\theta}-1\}$.

We now have 
$$Q^{(\widetilde{\theta})}(z)=\frac{B^{[\widetilde{\theta}]}(z)}{(1-z)^{\widetilde{\theta}}D(z)}.$$ 
By Lemma~\ref{lem:gamma}, we notice that \begin{align}\label{Dsym}
D(z)=(1-z)\sum_{n=0}^{\infty}\gamma_n(\0)z^n\sim {\sf Cst}\cdot (1-z)^{-\theta}\quad\text{ as $z\to 1^-$}.\end{align} 
We also note that as $z\to 1^-$,
\begin{align}\label{Bthetasym}B^{[\widetilde{\theta}]}(z)& ={
R^{[\widetilde{\theta}]}_B(z)= Q^{(\widetilde{\theta}-1)}(1^-)R^{[\widetilde{\theta}-1]}_D(z)- R_B^{[\widetilde{\theta}-1]}(z) \sim \left\{\begin{matrix} \textsf{Cst}\cdot\ln\big(\frac{1}{1-z}\big) & \text{if } \theta\in\N,\\ \textsf{Cst} &\text{otherwise.} \end{matrix}\right.}
\end{align}
Hence \eqref{sum.q2} is verified.

{Assume now that $\theta\in\N$. One can easily show by induction that
 $\frac{\text{d}}{\text{d}z} R_{B}^{[k]}(z)\sim \textsf{Cst}\cdot (1-z)^{-1}$ as $z\to 1^{-}$ for all $k\in\{0,1,\dots, \theta\}$.
Hence \begin{align}\label{DBsym}
\frac{\text{d}}{\text{d}z}B^{[\theta]}(z)=\frac{\text{d}}{\text{d}z}R_B^{[\theta]}(z)\sim \textsf{Cst}\cdot (1-z)^{-1} \quad\text{as $z\to 1^-$}.\end{align}
On the other hand, using \eqref{D.expand}, we have
\begin{align*}\frac{\text{d}}{\text{d}z}\big( (1-z)^{\theta}D(z) \big)&= (1-z)^{\theta}\left(-\theta\sum_{n=0}^{\infty}\gamma_n(\0)z^n+\sum_{n=0}^{\infty} n \big(1-g(f^{n}(\0)) \big)\gamma_n(\0) z^{n-1} \right).\end{align*}
By reason of \eqref{fn0} and \eqref{gf0}, we note that $1-g(f^n(\0))\sim \theta n^{-1}+O(n^{-(1+\eps)})$  for some $\eps\in (0,1)$  as $n\to\infty$. Hence
\begin{align}\label{DDsym}\frac{\text{d}}{\text{d}z}\big( (1-z)^{\theta}D(z) \big)=O\Big((1-z)^{\theta}\sum_{n=1}^{\infty}n^{-\eps}\gamma_n(\0) z^{n}\Big)=O((1-z)^{-(1-\eps)})  \quad\text{as } z\to 1^{-}\end{align}
where in the last equality we use Hardy–Littlewood tauberian theorem and the fact that $\gamma_n(\0)\sim C n^{\theta} $. Combining \eqref{Dsym}, \eqref{Bthetasym}, \eqref{DBsym} and \eqref{DDsym}, we have
$$\frac{\text{d}}{\text{d}z} Q^{(\theta)}(z)= \frac{(1-z)^{\theta}D(z)\frac{\text d}{\text{d}z}B^{[\theta]}(z) -B^{[\theta]}(z)\frac{\text{d}}{\text{d}z}\big( (1-z)^{\theta}D(z) \big) }{\big((1-z)^{\theta}D(z)\big)^2}\sim \textsf{Cst}\cdot (1-z)^{-1}\quad\text{as }z\to 1^-.$$ Hence \eqref{dQtheta} is verified.}\\

\underline{Part b}. 
Recall from the proof of Lemma~\ref{lem:Q} that $$F(s,n)= F(f(s),n-1)g(f(s))-F(\0,n-1)(g(f(s))-1) -\sum_{ r\in \mathcal{S}(N)}\P(U(n-1)=r) H_r(f(s)).$$
Differentiating both sides of the above equation at $s=\mathbf 1$, we obtain that 
$$\E [U(n)]= \Big( \E[U(n-1)]   -\mu_{n} N+\sum_{r\in \mathcal{S}(N)}\P(U(n-1)=r) (N-r)  \Big)\cdot\E[\psi].$$
{Iterating the above equality and taking the scalar product with $u$, we have}
\begin{equation}\label{eq.EU1}\E [\langle U(n),u\rangle ]=\langle K, u \rangle+\sum_{k=1}^{n-1} \Big(-\mu_{k} \langle N, u \rangle +\sum_{r\in \mathcal{S}(N)}\P(U(k-1)=r) \langle N-r, u \rangle\Big)\end{equation}
where we recall that $U(0)= K$ a.s. and $u$ is the right eigenvector of $\E[\psi]$ w.r.t. the maximal eigenvalue 1. 
{Taking $n\to\infty$ in \eqref{eq.EU1} and using \eqref{eq.Pi}}, we get
$$\sum_{n=0}^{\infty}\mu_{n}=Q(1^-) =\frac{\langle K,u\rangle+\sum_{r\in \mathcal{S}(N)}\Pi_r(1^-)\langle N-r,u\rangle}{\langle N,u\rangle}.$$
Hence we obtain
\begin{align*}
\E[\langle U(n), u\rangle]&= \langle N, u\rangle \sum_{k=n}^{\infty} \mu_{k}  -\sum_{k=n}^{\infty} \sum_{r\in \mathcal{S}(N)}\P(U(k-1)=r) \langle N-r, u \rangle
\\&\ge \langle N, u\rangle \Big[ \sum_{k=n}^{\infty} \mu_{k}  -\sum_{k=n}^{\infty} \sum_{r\in \mathcal{S}(N)}\P(U(k-1)=r)  \Big].
\end{align*}
By reason of \eqref{inq.pi1}, $\sum_{k=n}^{\infty} \sum_{r\in \mathcal{S}(N)}\P(U(k-1)=r)\le \mu_{n-1}=O(n^{-1-\theta})$. 
On the other hand, we notice that $\sum_{k=n}^{\infty} \mu_{k}\sim \frac{\varrho}{\theta} n^{-\theta}$ and $\langle U(n), u\rangle \le |U(n)|\cdot|u|=|U(n)|$. 
Hence $$\liminf_{n\to\infty}\frac{\E[| U(n)| \ |U(n)\neq \0]}{n}= \liminf_{n\to\infty} \frac{\E[|U(n)|]}{n\mu_{n}} \ge \liminf_{n\to\infty}  \frac{\langle N, u\rangle \Big(\frac{\varrho}{\theta} n^{-\theta}+O(n^{-1-\theta}) \Big) }{n\cdot\varrho n^{-\theta-1}}=\beta.$$
\end{proof}

\section{Phase transition for the speed of ERW}\label{sec:speed}
{We first prove the following result which in particular implies Part (a) of Theorem~\ref{th:main}.}
\begin{tcolorbox}
\begin{theorem}\label{thm.exist} Under Assumption \ref{Asum:A}, we have that
\begin{equation}\label{speed}
\lim_{n\to\infty}\frac{X_n}{n}= \frac{1}{1+ \langle\E[Z_{\infty}], \varsigma\rangle }\quad\text{almost surely},
\end{equation}
where $\varsigma=(2,1,\dots,1)$.
\end{theorem}
\end{tcolorbox}
\begin{proof}
{
In virtue of Proposition \ref{ergodic}, the distribution of $Z_{\infty}$ does not depend on $Z_0$ and
$$\lim_{n\to\infty} \sum_{i=0}^{n-1}\frac{\langle Z_i, \varsigma\rangle}{n}= \left\langle\E[Z_{\infty}], \varsigma\right\rangle \text{ almost surely}.$$
Assume w.l.o.g. that $Z_0=0$ almost surely. By reason of Proposition \ref{Markovian},  $\sum_{i=0}^{n-1}\langle Z_i, \varsigma\rangle$  is equal to $\sum_{i=0}^{n-1}\langle V^n_i, \varsigma\rangle$ in distribution. Moreover, recall from \eqref{hitting} that $T_n\sim n+\sum_{i=0}^{n-1} \langle V^n_i, \varsigma\rangle$ as $n\to\infty$. Hence 
$$\lim_{n\to\infty} \frac{T_n}{n}= 1+\langle\E[Z_{\infty}], \varsigma\rangle\quad \text{in probability}.$$}
{Using an argument by Zerner (see the proof of Theorem 13 in \cite{Z2005}), we next show that $X_n/n$ and $\lim_{n\to\infty}n/T_n$ converge in probability to the same limit. Indeed, for $n\ge 0$, set
$S_n=\sup\{k: T_k\le n\}$. We have 
$T_{S_n}\le n< T_{S_n+1}$. It immediately follows that 
\begin{equation} \label{eq:samelim}\lim_{n\to\infty}\frac{n}{S_n}=\lim_{n\to\infty} \frac{T_{S_n}}{S_n}=\lim_{n\to\infty} \frac{T_n}{n}=1+\langle\E[Z_{\infty}], \varsigma\rangle \quad\text{in probability}.\end{equation}
Since $n< T_{S_n+1}$ and before time $T_{S_n+1}$, the walk is always below level $S_n+1$, we must have $X_n\le S_n$. On the other hand
$X_{n}=S_n+X_n-X_{T_{S_n}}\ge S_n -L(n-T_{S_n}).$ It follows that 
\begin{align*}\label{lim2}\frac{S_n}n - L\left(1- \frac{T_{S_n}}{n}\right) \le\frac{X_n}n \le \frac{S_n}{n}.\end{align*}
Using \eqref{eq:samelim}, we note that $$ \lim_{n\to\infty}\frac{T_{S_n}}{n} =\lim_{n\to\infty} \frac{T_{S_n}}{S_n} \frac{S_n}n = 1\quad \text{in probability}.$$
 As a result, $$\lim_{n\to\infty}\frac{X_n}{n}=\lim_{n\to\infty}\frac{S_n}{n}=\lim_{n\to\infty}\frac{n}{T_n}= \frac{1}{1+\langle\E[Z_{\infty}], \varsigma\rangle} \quad \text{in probability}.$$
Furthermore, $X_n/n$ converges almost surely as $n\to\infty$ thanks to Proposition 3.1 in \cite{CHN2021}. We note that the proof of Proposition 3.1 in \cite{CHN2021} generally holds for transient multi-excited random walk with long jumps both to the left and the right assuming that $$\sup\Big(   \operatorname \supp(\nu)\cup \bigcup_{i=2}^M \supp(q_i)\Big)\le \sup( \supp(q_1))<\infty$$
where $\supp(\mu)$ stands for the support of measure $\mu$.
This ends the proof of the theorem.}
\end{proof}

Theorem \ref{thm.exist} also implies that Part (b) of Theorem~\ref{th:main} is equivalent to the following theorem.
\begin{tcolorbox}
\begin{theorem}\label{thm:speed} 
$\E[Z_{\infty,\l}]<\infty$ for all $\l \in [L]$ if and only if $\delta>2.$ 
\end{theorem}
\end{tcolorbox}

To prove the above theorem, we will need some preliminary results. The next proposition is immediate from the proof of Proposition 3.6 in \cite{DS2008}.
\begin{proposition}\label{lem:genG}
Suppose that for $s\in [0,1]$
\begin{equation}\label{feq.G}1-G\left(\frac{1}{2-s}\right)=a(s)(1-G(s))+b(s)\end{equation}
where 
\begin{itemize}
\item[I.] $a(s)$ and $b(s)$ are analytic functions in some neighborhood of 1 such that $a(1)=1, a'(1)=\delta$ for some $\delta>1$ and $b(1)=0$;
\item[II.] $G$ is a function defined on $[0,1]$ such that $G$ is left continuous function at $1$, $G'(1^-)\in (0,\infty]$ and there exists $\epsilon\in (0,1)$ such that $G^{(i)}(s)>0$ for each $s\in (1-\epsilon,1)$ and $i\in \N$.
\end{itemize}
   Then, the following statements hold true:
\begin{itemize}
\item[(i)] $b'(1)=0$.
\item[(ii)] If $\delta>2$ then $b''(1)>0$ and $1-G(1-s)=\frac{b''(1)}{2(\delta-2)}s+O(s^{2\wedge (\delta-1)}) $ as $s\downarrow 0$.
\item[(iii)] If $\delta=2$ and $b''(1)=0$ then
$G^{(i)}(1^{-})<\infty$ for all $i\in \mathbb{N}$.
\item[(iv)] If $\delta=2$ and $b''(1)\neq 0$ then $1-G(1-s)\sim C s |\ln(s)|$ as $s\downarrow 0$ for some constant $C>0$.
\end{itemize}
\end{proposition}

In this section, we consider the following function \begin{align}\label{def.G}G(s):=\E\Big[ \prod_{\l=1}^L (1+\l(s-1))^{Z_{\infty,\l}} \Big], \ s\in[0,1].
\end{align}
Notice that $G'(1)=\E\Big[ \sum_{\l=1}^L \l Z_{\,\infty, \l}\Big] $. For $s_0\in \left(\frac{L-1}{L},1\right]$, we have
\begin{align*}G(s)&=\E\Big[ \prod_{\l=1}^L \left(1-\l(1-s_0)+\l(s-s_0)\right)^{Z_{\infty,\l}} \Big]\\& =\sum_{k\in \Z_+^L}\P(Z_{\infty}=k)\prod_{\l=1}^L\sum_{j=0}^{k_{\l}}\binom{k_{\l}}{j}(1-\l(1-s_0))^j (\l(s-s_0))^{k_{\l}-j}.
\end{align*}
It follows that $G(s)$ can be expanded as a power series of $s-s_0$ with all positive coefficients. As a consequence, $G^{(i)}(s)>0$ for all $s\in \left(\frac{L-1}{L},1\right)$ and $i\in \N$. 
Hence the condition II of Proposition~\ref{lem:genG} is verified. We next show (in Proposition~\ref{prop:1} and Proposition \ref{prop:2} below) that there exist functions $a$ and $b$ such that the condition I and the functional equation \eqref{feq.G} (with $G$ given by \eqref{def.G}) are fulfilled.

Recall that $(\eta_n)_{n\ge1}$ is a sequence of i.i.d. random vectors with multivariate geometrical law defined in \eqref{geom}. For $\l \in [L]$, set $\rho_{\l}=\nu(-\l)/\nu(1)$. Notice that the probability generating function of $\eta_1$ is given by
\begin{align}\label{characEta}\E\Big[ \prod_{\l=1}^L s_{\l}^{\eta_{1,\l}}\Big]=\frac{1}{1+\sum_{\l=1}^L\rho_{\l}(1-s_{\l})} \end{align}
{for $s=(s_1,s_2,\dots, s_L)$ such that $\sum_{\ell=1}^Ls_{\ell}\nu(-\ell)<1$.}

\begin{lemma}\label{pr:AM} For $\l\in[L]$,
 \begin{align*}
 \E\left[A_{\ell}({M-1}) \right]=\sum_{i=1}^M \left( q_{i}(-\l)+\rho_{\l} (1-q_i(1))\right).
 \end{align*}
 \end{lemma}
\begin{proof}
Set $S=\sum_{i=1}^M \xi_{i,L+1}$. We have that $\E[S]=\sum_{i=1}^Mq_i(1),$
and
$$A_{\ell}({M-1})= \sum_{i=1}^M \xi_{i,\ell}+ \sum_{i=M+1}^{\gamma_{M-1} } \xi_{i,\l},$$
On the other hand, recall that 
$$\gamma_{M-1}=\inf\left\{k\ge 1: \sum_{i=1}^k \xi_{i,L+1}=M\right\}=M+{\inf \left\{k\ge 1: \sum_{i=M+1}^{M+k} \xi_{i,L+1}=M-S\right\}.}$$
Hence $\sum_{i=M+1}^{\gamma_{M-1} } \xi_{i,\l}$ has the same distribution as
$\sum_{i=1}^{M-S} \eta_{i,\l}.$ 
Therefore
$$\E[A_{\ell}({M-1})]=\sum_{i=1}^M\E[ \xi_{i,\l}] +\E[\eta_{1,\l}]\E[M-S] =\sum_{i=1}^M q_{i}(-\l)+\rho_{\l} \sum_{i=1}^M(1-q_i(1)). $$
\end{proof}
{We can combine Proposition \ref{iden.A} with Proposition~\ref{pr:AM} to compute {$\E[A(m)]$} for $m \ge M$. We also note from \eqref{characEta} that $\E[\eta_{1,\ell}]=\rho_{\ell}$ for $\ell\in [L]$. As a result, for $\l\in[L ]$ and $k\in \Z^L_+$ such that $|k|=k_1+\dots +k_L\ge M-1$, we have
\begin{align*}\E[Z_{1,\l}|Z_{0}=k]&  =\E\left[A_{\l}(|k|)\right]+k_{\l+1}  = \E\left[A_{\l}(M-1)\right]+\sum_{i=1}^{|k|-M+1}\E[\eta_i]+k_{\l+1}\\
&=
 \sum_{i=1}^M q_{i}(-\l)+\rho_{\l}\left(|k|+1- \sum_{i=1}^Mq_i(1)\right)+ k_{\l+1}
\end{align*}
where we use the convention $k_{L+1}=0$.}

\begin{proposition}\label{prop:1}The function $G$ defined by \eqref{def.G}  satisfies the functional equation \eqref{feq.G}
where we define
\begin{align*}a(s)&:= \frac{1}{\E\left[ \prod_{\l=1}^L(1+\l(s-1))^{A_{\l}(M-1)}\right](2-s)^{M-1}},\\
 b(s)&:=\sum_{|k|\le M-2}\P(Z_{\infty}=k)\left(a(s)\E\Big[ \prod_{\l=1}^L (1+\l(s-1))^{A_{\l}(|k|)+k_{\l+1}}\Big] -\frac{ \prod_{\l=1}^{L} (1+\l(s-1))^{k_{\l+1}}}{(2-s)^{|k|}}\right)\\
& -a(s)+1.
 \end{align*}
\end{proposition}
\begin{proof}
We have that
\begin{align*}
& G(s)=\sum_{k\in \Z^L_+}\P\left(Z_{\infty}=k \right)\E\Big[ \prod_{\l=1}^L (1+\l(s-1))^{Z_{1,\l}}\ \big|\ Z_0=k\Big].
\end{align*}
Recall that given $\{Z_0=(k_1,k_2,\dots, k_L)\}$, the random vector $Z_1=(Z_{1,1},Z_{1,2},\dots, Z_{1,L})$ has the same distribution as $\left( A_1(|k|)+k_2,\dots , A_{L-1}(|k|)+k_L ,  A_L(|k|)\right)$. Using Proposition \ref{iden.A}, we thus have
\begin{align*}
G(s)&= \sum_{|k|\le M-2 }\P\left(Z_{\infty}=k \right)\E\Big[ \prod_{\l=1}^L (1+\l(s-1))^{A_{\l}(|k|)}\Big]\prod_{\l=1}^{L-1} (1+\l(s-1))^{k_{\l+1}} \\
& +\sum_{|k|\ge M-1 }\P\left(Z_{\infty}=k \right)\E\Big[ \prod_{\l=1}^L(1+\l(s-1))^{A_{\l}(M-1)+\eta_{1,\l}+\dots+\eta_{|k|-M+1,\l}}\Big]\prod_{\l=1}^{L-1} (1+\l(s-1))^{k_{\l+1}}.
\end{align*}
On the other hand, using \eqref{characEta} and the fact that $\sum_{\l=1}^L\l \rho_{\l}=1$, we obtain
$$\E\Big[ \prod_{\l=1}^L (1+\l(s-1))^{\eta_{1,\l}}\Big]=\frac{1}{2-s}. $$
Hence
\begin{align*}& \sum_{k\in \Z^L_+ } \P\left(Z_{\infty}=k \right) \E\Big[ \prod_{\l=1}^L (1+\l(s-1))^{A_{\l}(M-1)+\eta_{1,\l}+\dots+\eta_{|k|-M+1,\l}}\Big] \prod_{\l=1}^{L-1} (1+\l(s-1))^{k_{\l+1}}
\\
& =  \sum_{k\in \Z^L_+ } \P\left(Z_{\infty}=k \right)  \frac{\E\Big[ \prod_{\l=1}^L(1+\l(s-1))^{A_{\l}(M-1)}\Big] \prod_{\l=1}^{L-1} (1+\l(s-1))^{k_{\l+1}}}{\left(2-s\right)^{|k|-M+1}} \\
& = \E\Big[ \prod_{\l=1}^L(1+\l(s-1))^{A_{\l}(M-1)}\Big](2-s)^{M-1}\sum_{k\in \Z^L_+ } \P\left(Z_{\infty}=k \right) \prod_{\l=1}^{L} \left(1+\l\left(\frac{1}{2-s}-1\right)\right)^{k_{\l}}\\
& =\frac{1}{a(s)} G\left(\frac{1}{2-s}
\right).
\end{align*} 
Therefore
\begin{align*}
G(s)&= \sum_{k\in \Z^L_+ } \P\left(Z_{\infty}=k \right)  \frac{\E\Big[ \prod_{\l=1}^L (1+\l(s-1))^{A_{\l}(M-1)}\Big] \prod_{\l=1}^{L-1} (1+\l(s-1))^{k_{\l+1}}}{\left( 2-s\right)^{|k|-M+1}} \\
& +\sum_{ |k|\le M-2 }\P\left(Z_{\infty}=k \right)\left(\E\Big[ \prod_{\l=1}^L (1+\l(s-1))^{A_{\l}(|k|)}\Big] -\frac{ \E\Big[ \prod_{\l=1}^L (1+\l(s-1))^{A_{\l}(M-1)}\Big] }{ \left( 2-s\right)^{|k|-M+1} }\right)\\
& \times \prod_{\l=1}^{L-1} (1+\l(s-1))^{k_{\l+1}}  = \frac{1}{a(s)} G\left(\frac{1}{2-s}
\right) + 1+\frac{b(s)-1}{a(s)}.
\end{align*}
\end{proof}
Recall from \eqref{eq:totdr} that the expected total drift $\delta$ of the cookie environment is given by
$$\delta:=\sum_{j=1}^M\left( q_j(1)-\sum_{\l=1}^L \l q_{j}(-\l) \right).$$
We must have that $\delta>1$ as $\X$ is assumed to be transient to the right.

\begin{proposition}\label{prop:2} The functions $a(s)$ and $b(s)$ defined in Proposition~\ref{prop:1} satisfies the condition (I) of Proposition~\ref{lem:genG}. More specifically,
\begin{align*}a(1-s) = 1- \left(\delta -1 \right) s +o(s)\quad \text{and}\quad b(1-s) =b'(1)s+o(s)\quad \text{ as } s\to 0
\end{align*}
where $$b'(1)=(\delta-1)- \sum_{|k|\le M-2}\P[Z_{\infty}=k]  \left( \delta-1-|k|+ \sum_{\l=1}^L \l \E[A_{\l}(|k|)] \right). $$
\end{proposition}
\begin{proof}
Using Taylor's expansion, we have
\begin{align*}
 a(1-s)&=\frac{1}{\E\left[ \prod_{\l=1}^L (1-\l s )^{A_{\l}(M-1)}\right] \left( 1+s \right)^{M-1}} = 1- \left( M-1- \sum_{\l=1}^L  \l \E[A_{\l}(M-1)]\right) s+o(s)\\
& = 1- \left(\sum_{i=1}^M \left( q_i(1)-\l q_{i}(-\l)\right) -1 \right) s +o(s)=1-(\delta-1)s+o(s)
\end{align*}
as $s\to 0$, {in which the third identity follows from Lemma \ref{pr:AM} and the fact that $\sum_{\ell=1}^L\ell\rho_{\ell}=1$}.
Furthermore,
\begin{align*}& b(1-s)=1-a(1-s)+\\
& +\sum_{|k|\le M-2}\P[Z_{\infty}=k]\left(a(1-s)\E\Big[ \prod_{\l=1}^L (1-\l s)^{A_{\l}(|k|)}\Big] - \left( 1+s\right)^{-|k|}\right)\prod_{\l=1}^{L-1} (1- \l s)^{k_{\l+1}}\\
&= (\delta-1)s - \sum_{|k|\le M-2}\P[Z_{\infty}=k]  \left( \delta-1-|k|+ \sum_{\l=1}^L \l\E[A_{\l}(|k|)] \right)s +o(s).
\end{align*}
\end{proof}
\begin{remark}
Note that the functional equation \eqref{feq.G} has \textbf{the same} form with the one in \cite{DS2008}. However the coefficient function $b$ defined in Proposition \ref{prop:1} is more complicated and strongly depends on the distribution of $Z_{\infty}$ while the function $G$ defined by \eqref{def.G} is \textbf{not} a probability generating function when $L\ge 2$ and it will not give us the full information to compute $b$. Nevertheless, Proposition~\ref{lem:genG} still play a crucial role in the proof of Theorem~\ref{th:main}. 
\end{remark}

Let $\eta=(\eta^1,\eta^2,\dots, \eta^{L})$ be a $L$-dimensional random vectors with multivariate geometric law defined by
\begin{align*}
\P\left(\eta=(i_1,i_2,\dots,i_{L})\right)= \frac{\nu(1) }{(i_1+i_2+\dots+i_L)!} \prod_{k \in [L]}  i_k ! \nu(-k)^{i_{k}},
\end{align*}
for each $i=(i_1,i_2,\dots, i_L)\in \Z^L_+$. Recall that the probability generating function of $\eta$ is given by
\begin{align*}\E\Big[ \prod_{\l=1}^L s_{\l}^{\eta^{\l}}\Big]=\frac{1}{1+\sum_{\l=1}^L\rho_{\l}(1-s_{\l})} \quad \text{with }\rho_{\l}=\nu(-\l)/\nu(1).\end{align*}
In particular, we have $\E[\eta]=(\rho_1,\rho_2,\dots, \rho_{L})$. 

Let $(\vartheta(k,n))_{k,n\in \N}=(\vartheta_{i,j}(k,n), i,j\in [L] )_{k,n\in \N}$ be a sequence of i.i.d. $L\times L$ random matrices such that its rows are i.i.d. copies of $\eta$. Let $(W(n))_{n\ge 0}=(W_1(n),\dots, W_L(n))_{n\ge0}$ be a multi-type branching process with $(N_1,N_2,\dots, N_L)$-emigration such that for each $n\in\N$ and $j\in[L]$,
\begin{equation}\label{mult.Y}
W_{j}(n)=\sum_{i=1}^{L}\sum_{k=1}^{\varphi_i(W(n-1))}\chi_{i,j}(k,n)
\end{equation}
where $$\chi_{i,j}(k,n)=\vartheta_{i,j}(k,n) +\delta_{i-1,j}$$ (here $\delta_{i,j}$ stands for the Kronecker delta) and {$\varphi_i(w):=(w_i-N_i)\mathds{1}_{\{w_{\l} \ge N_{\ell},\ \forall \l \in[L] \}}$}.
\begin{lemma}\label{lem:critical.Y}  Assume $\nu(-L)>0$. Then
 $(W(n))_{n\ge 0}$ is a critical multi-type Galton-Watson process with $(N_1, N_2,\dots, N_L)$-emigration according to Definition \ref{def:emigration} and Assumption \ref{Asum:B}. 
\end{lemma}
\begin{proof}
We have $$\overline{\chi}:= \E\left[(\chi_{i,j}(1,1))_{i,j\in[L]}\right]=\begin{pmatrix}
 \rho_1 & \rho_2 & \dots & \rho_{L-1} & \rho_L \\
 \rho_1+1 & \rho_2 & \dots & \rho_{L-1} & \rho_L \\
 \rho_1 & \rho_2+1 & \dots  & \rho_{L-1} & \rho_L \\
 \vdots & \vdots & \ddots  & \ddots & \vdots \\
 \rho_1 & \rho_2 & \dots &   \rho_{L-1}+1 & \rho_L \end{pmatrix}$$
 and notice that all the entries of
 $$\begin{pmatrix}
 0 & 0 & \dots & 0 & \rho_L \\
 1 & 0 & \dots & 0 & \rho_L \\
 0 & 1 & \dots  & 0 & \rho_L \\
 0 & \vdots & \ddots  & \ddots & \vdots \\
 0 & 0 & \dots &  1 & \rho_L \end{pmatrix}^L$$
are positive as $\rho_L=\nu(-L)/\nu(1)>0$. 
Hence $\overline{\chi}$ is positively regular.  
 
Applying the determinant formula {(see, e.g., Theorem 18.1.1 in \cite{Harville})} $$ \det \left(\Sigma +\mathbf{x}.\mathbf{y} ^{\textsf{T}}\right)=\left(1+\mathbf{x}^{\textsf{T}} \Sigma^{-1}\mathbf{y} \right)\det(\Sigma)$$ (in which $\Sigma$ is an invertible matrix, $\mathbf x$ and $\mathbf y$ are column vectors) to $\mathbf x=(\rho_1,\rho_2,\dots,\rho_L)^{\textsf {T}}$, $\mathbf y=(1,1,\dots, 1)^{\textsf {T}}$,
$$\Sigma=\begin{pmatrix}
 -\lambda & 0       & \dots & 0 & 0 \\
     1    & -\lambda &\dots & 0  & 0\\
     0    &  1       & \dots & 0 & 0\\
 \vdots & \vdots & \ddots  & \ddots & \vdots\\
     0    &  0       & \dots &  1 & -\lambda
 \end{pmatrix}\quad\text{and} \quad \Sigma^{-1}=\begin{pmatrix}
-\frac{1}{\lambda} & 0 & \dots & 0 \\
 -\frac{1}{\lambda^2} & -\frac{1}{\lambda} & \dots & 0 \\
 \vdots & \vdots & \ddots & \vdots \\
 -\frac{1}{\lambda^L} & -\frac{1}{\lambda^{L-1}} & \dots & -\frac{1}{\lambda} 
\end{pmatrix},$$
 we obtain
 {$$\Phi(\lambda):=\det(\overline{\chi}-\lambda I)=(-1)^L\left(\lambda^{L}- \sum_{j=0}^{L-1}\left(\sum^{L}_{\l=L-j}\rho_{\l}\right)\lambda^{j}\right)=(-1)^L \sum_{j=0}^{L-1}\left(\sum^{L}_{\l=L-j}\rho_{\l}\right)(\lambda^{L}-\lambda^{j})$$
where the last equality follows from the fact that $$\sum_{j=0}^{L-1}\sum^{L}_{\l=L-j}\rho_{\l}=\sum_{\l=1}^L\l \rho_{\l}=1.$$ 
As a result, $\Phi(1)=0$ and \begin{align*}
  \Phi'(1)= (-1)^L \sum_{j=0}^{L-1}(L-j)\left(\sum^{L}_{\l=L-j}\rho_{\l}\right)\neq 0.
\end{align*}
yielding that $\lambda=1$ is a simple eigenvalue of $\overline{\chi}$. Using the triangle inequality, we have that $$|\Phi(\lambda)|\ge|\lambda^{L}|-\sum_{j=1}^{L-1}\left(\sum^{L}_{\l=L-j}\rho_{\l}\right)|\lambda^j|=  \sum_{j=1}^{L-1}\left(\sum^{L}_{\l=L-j}\rho_{\l}\right)(|\lambda|^{L}-|\lambda|^j)>0 \quad \text{if} \quad |\lambda|>1.$$
Therefore all the eigenvalues of $\overline{\chi}$ must lie in the closed unit disk.} Hence $\lambda=1$ is also the largest eigenvalue in modulus. This ends the proof of the lemma.
\end{proof}
 
\begin{remark}
The right and left eigenvectors of the maximal eigenvalue $\lambda=1$ are given respectively by
 \begin{align*}u&=\frac{2}{L(L+1)}\left(1,2,\dots, L-1,L\right),\\
 v & =\frac{1}{\rho_L+\frac{1}{L+1}\sum_{\l=1}^{L-1}\l(\l+1){\rho_{\l}}} \left({\rho_L}+L\sum_{\l=1}^{L-1} \rho_{L-\l},\dots, {\rho_L}+{L(\rho_{L-2}+\rho_{L-1})},{\rho_L}+{L\rho_{L-1}},{\rho_L}\right).
 \end{align*}
It is also clear that \eqref{cond:moment} holds true since $\E[\chi_{ij}^k]<\infty $ for all $k\ge 1$ and $i,j\in [L]$. 
\end{remark}

\begin{proposition}\label{lem:moment}
Assume that $\delta=2$. Then there exists a positive integer $\kappa$ such that
$$\E[|Z_{\infty}|^{\kappa}]=\infty.$$
\end{proposition}
\begin{proof}
Define $\tau=\inf\{n\ge1 :\ Z_n=\0\}$. Notice that $\E[\tau| Z_0=\0]<\infty$ as $\ZZ$ is positive recurrent. Furthermore, for any function $\pi:\Z^L\to \R_+$, we have (see e.g. Theorem 1.7.5 in \cite{Norris})
$$\E[\pi(Z_{\infty})]={\frac{\E\Big[\sum_{n=0}^{\tau-1}\pi(Z_n)\ | Z_0=\0\Big]}{\E[\tau | Z_0=\0]}}.$$ 
{Let $L'=\max\{\ell\in[L]:\nu(-\ell)>0\}$ and $z^*=(z^*_1,\dots,z^*_L)$ with $z^*_{\l}=M$ for $1\le \l\le L'$ and $z^*=0$ for $L'+1\le \l\le L$.} Let $\kappa$ be a fixed positive integer that we will choose later. By setting $\pi(z)=\left(\sum_{\ell=1}^L  z_{\ell}\right)^{\kappa}$, we obtain
\begin{align}\nonumber \E[|Z_{\infty}|^{\kappa}]& {= \frac{1}{\E[\tau | Z_0=\0]}\E\Big[\sum_{n=0}^{\tau-1}|Z_n|^{\kappa}\  | Z_0=\0 \Big]}\\
\label{ineq1} & {\ge \frac{\P\big(Z_L=z^*, \tau>L | Z_0=\0\big)}{\E[\tau | Z_0=\0]}\E\Big[\sum_{n=0}^{\infty}|Z_{n\wedge \tau}|^{\kappa}\ | Z_0=z^*\Big].}
\end{align}
Using \eqref{Zn}, we note that $\P\big(Z_L=z^*| \tau>L , Z_0=\0\big)=\P\big(Z_L=z^*| Z_k\neq \0 \text{ for all }1\le k\le L-1 , Z_0=\0\big)>0$ and thus \begin{align}\label{Pzstar}
    \P\big(Z_L=z^*, \tau>L | Z_0=\0\big)>0.
    \end{align}

We next use a coupling argument to estimate the order of $\E[|Z_{n\wedge \tau}|^{\kappa}| Z_0=z_*]$ as $n\to\infty$. {Recall from Remark~\ref{rem:branching} that $$Z_{n}=\left\{\begin{matrix}\displaystyle A^{(n)}(ML-1)+\sum_{k=(L-1)M+1}^{|Z_{n-1}|-M+1}\eta^{(n)}_{k} + \widetilde{Z}_{n-1} & \text{if}& |Z_{n}|\ge ML-1,\\ 
A^{(n)}(|Z_{n-1}|)+ \widetilde{Z}_{n-1} & \text{if} & |Z_{n}|<ML-1\end{matrix}\right.$$
 where we denote $\widetilde{z}=(z_2,z_3,\dots, z_{L+1},0)$ for each $z=(z_1,z_2,\dots, z_L)\in \R^L$; and $\eta_k^{(n)}=(\eta_{k,1}^{(n)},\dots, \eta_{k,L}^{(n)})$ with $k, n\in \N$ are i.i.d. copies of the random vector $\eta$ which are independent of $Z_0$. Let $(W(n))_{n\ge0}$ be the multi-type branching process with $N$-emigration defined by \eqref{mult.Y}, in which we set $N=(N_1,N_2,\dots,N_L)=(M-1,M,\dots,M)$, $W(0)=(M,M,\dots, M)$ and $\vartheta_{i,j}(k,n)=\eta_{(L-1)M+k,j}^{(n)}$ for $i,j\in[L]$ and $k, n\in\N$.} For the sake of simplicity, we assume from now on that $L=L'$ (otherwise, by reducing the dimension of $W(n)$, one can easily handle the case $L'<L$ by the same argument used in the case $L=L'$). In this case, we note that $z^*=(M,M,\dots,M)$. 

Recall that for $L$-dimensional vectors $x$ and $y$, we write $x\succeq y$ if ${x}_{\l}\ge y_{\l}$ for all $\l\in [L]$. Assuming $Z_0=z^*=(M,M,\dots,M)$, we will show that \begin{equation}\label{sto.dom}Z_n {\succeq} {W}(n)\quad\text{for all  } n\ge 1.\end{equation}
Indeed, for $n=1$, we have {$$Z_1=  A^{(1)}({LM-1}) + \eta^{(1)}_{(L-1)M+1}+(M,M,\cdots, M,0)\ge \eta^{(1)}_{(L-1)M+1}=W(1).$$ Suppose that $Z_{n-1}{\succeq} W(n-1)$ for some $n\in\N$. If $Z_{n-1,\l}< N_{\l}$ for some $\l \in[L]$ then $W_{\l}(n-1)\le Z_{n-1,\ell} < N_{\ell}$ and thus $W(n)=\0$. On the other hand, if $Z_{n-1,\l}\ge N_{\l}$ for all $\l\in [L]$ then $|Z_{n-1}|\ge ML-1$ and thus
 for $\ell\in[L]$,\begin{align*}Z_{n,\ell}&\ge \sum_{k=(L-1)M+1}^{|Z_{n-1}|-M+1}\eta^{(n)}_{k,\ell} + {Z}_{n-1,\ell+1}-N_{\ell+1}= \sum_{i=1}^L\sum_{k=1}^{Z_{n-1,i}-N_i}\left(\vartheta_{i,\ell}(k,n) + \delta_{i-1,\ell} \right) \\
 &  =\sum_{i=1}^L\sum_{k=1}^{Z_{n-1,i}-N_i}\chi_{i,\ell}(k,n)\ge W_{\ell}(n),
 \end{align*}
 in which we use the convention that $Z_{n-1,L+1}=N_{L+1}=0$.
Hence $Z_n \succeq W(n)$. By the principle of mathematical induction, we deduce \eqref{sto.dom}. It implies additionally that on the event $\{Z_0=z^*\}$, $W(n)=\0$ for all $n\ge \tau$.
Hence, for all $n\in\N$,
\begin{align}\label{ineq2}\E\left[|Z_{n\wedge \tau}|^{\kappa}\ | Z_0=z^*\right] \ge \E\big[\big|{W}({n})\big|^{\kappa}\ | Z_0=z^*\big]=\E[|W(n)|^{\kappa}]
\end{align}
where the last equality follows from the fact that  $(W(n))_{n\ge0}$ is independent of $Z_0$.}

On the other hand, $(W(n))_{n\ge0}$ is a critical multi-type Galton-Watson process with $(M-1,M,\dots, M)$-emigration. By Theorem~\ref{lem:branching}, there exist positive constants $c_1, c_2$ and $\theta$ such that 
 $$\P\left({W}(n)\neq \0\right)\ge \frac{c_1}{n^{\theta+1}}, \quad \E\left[\big|{W}(n)\big|\ |\ {W}(n)\neq \0\right]\ge {c_2}{n}.$$
 Choose $\kappa=\left\lfloor\theta\right\rfloor+1$. Using Jensen inequality, we thus have \begin{align}
\nonumber \E\left[\big|{W}(n)\big|^{\kappa}\right]& =\E\left[\big|{W}(n)\big|^{\kappa}\ |\ {W}(n) \neq \0\right]\P\left({W}(n)\neq \0\right)\\
 & \ge \E\left[\big|{W}(n)\big|\ | \  W(n) \neq \0\right]^{\kappa}\P\left({W}(n)\neq \0\right) 
\label{ineq3} \ge \frac{c_2^{\kappa} c_1 }{n^{\theta-\left\lfloor\theta\right\rfloor}}.
 \end{align}
 Combining \eqref{ineq1}, \eqref{Pzstar}, \eqref{ineq2} and \eqref{ineq3}, we obtain that $\E[|Z_{\infty}|^{\kappa}]=\infty.$
\end{proof}
We now turn to the proof of our main result.
\begin{proof}[Proof of Part (b), Theorem~\ref{th:main}]
Remind that this part is equivalent to Theorem~\ref{thm:speed}.
As the function $G$ defined by \eqref{def.G} satisfies the functional equation \eqref{feq.G} and the conditions I, II of Proposition~\ref{lem:genG},  
it follows  from Proposition \ref{lem:genG}(ii) that if $\delta>2$ then $$G'(1^-)=\E\Big[\sum_{\l=1}^L \l Z_{\infty,\l}\Big]=\frac{b''(1)}{2(\delta-2)}<\infty.$$ The above fact and \eqref{speed} imply that the random walk $\X$ has positive speed in the supercritical case $\delta>2$. 

Let us now consider the critical case $\delta=2$. If $b''(1)\neq 0$ then by the virtue of Proposition \ref{lem:genG}(iv), we must have $G'(1^-)=\infty$. Hence, to prove that $G'(1^-)=\infty$, it is sufficient to exclude the case $b''(1)=0$. Assume now that $b''(1)=0$. By Proposition \ref{lem:genG}(iii), $G^{(i)}(1^-)<\infty$ for all $i\in \N$. On the other hand, by Proposition \ref{lem:moment}, there exists a positive integer $\kappa$ such that $\E[|Z_{\infty}|^{\kappa}]=\infty$ and thus $G^{(\kappa)}(1^-)=\infty$, which is a contradiction. Hence $G'(1^-)=\infty$ and thus there exists $\ell\in[L]$ such that $\E[Z_{\infty,\ell}]=\infty$. It follows that a.s. $\lim_{n\to\infty }X_n/n=0.$  

The subcritical case can be solved by showing the monotonicity of the speed as follows. Assume that $1<\delta<2$. There exist probability measures $\widehat{q}_1, \widehat{q}_2,\dots, \widehat{q}_M$ on $\{-L,-L+1,\dots,-1,1\}$ such that 
$\widehat q_j(-\ell)\le q_j(-\ell)$ for each $\ell\in [L]$, $j\in[M]$
and
$$\widehat \delta:=\sum_{j=1}^M\left( \widehat q_j(1)-\sum_{\l=1}^L \l \widehat q_{j}(-\l) \right)=2.$$
Let $\widehat{\X}=(\widehat X_n)_n$ be the $(L,1)$-excited random walk w.r.t the cookie environment $\widehat \omega$ defined by
$$\widehat \omega(j,i) =\left\{\begin{array}{ll} \widehat q_{j}(i), \quad & \text{if}\ 1\le j \le M,\\
\nu(i), \quad  & \text{if}\ j >  M. \end{array}\right.$$
We thus have that a.s. $\lim_{n\to\infty}\widehat X_n/n=0$. Let $\mathbf Z$ and $\widehat{\mathbf Z}$ be respectively the Markov chains associated with $\X$ and $\widehat \X$ as defined by \eqref{def.Z}. Let $Z_{\infty}$ and $\widehat Z_{\infty}$ be their limiting distributions. For $h, \widehat{h}, k\in \Z_+^L$ with $ h \succeq \widehat{h}$, we notice that
$$\P( \widehat{Z}_n\succeq k |\widehat Z_{n-1}=\widehat{h} )\le \P( \widehat{Z}_n\succeq k |\widehat Z_{n-1}=h ) \le \P( Z_n\succeq k |Z_{n-1}=h ).$$
Applying Strassen's theorem on stochastic dominance for Markov chains (see e.g. Theorem 5.8, Chapter IV, p. 134 in \cite{Lindvall} or Theorem 7.15 in \cite{Holander}), we have that $\widehat{\mathbf Z}$ is stochastically dominated by $\mathbf Z$. In particular, $\E[\widehat{Z}_{\infty}]\le \E[{Z}_{\infty}]$. 
Combining the above fact and the speed formula \eqref{speed}, we conclude that a.s. $$\lim_{n\to\infty}\frac{ X_n}{n}\le \lim_{n\to\infty}\frac{\widehat X_n}{n}=0.$$
This ends the proof of our main theorem. 
\end{proof}

\section*{Acknowledgement}
The author would like to thank Andrea Collevecchio, Kais Hamza and the anonymous referees for their thorough reading and their constructive suggestions which improved the quality of the manuscript.

\bibliographystyle{amsplain}
\bibliography{mybib}

\providecommand{\bysame}{\leavevmode\hbox to3em{\hrulefill}\thinspace}
\providecommand{\MR}{\relax\ifhmode\unskip\space\fi MR }
\providecommand{\MRhref}[2]{%
  \href{http://www.ams.org/mathscinet-getitem?mr=#1}{#2}
}
\providecommand{\href}[2]{#2}
\begin{thebibliography}{10}

\bibitem{AN1972}
K.~B. Athreya and P.~E. Ney, \emph{Branching processes}, Die Grundlehren der
  mathematischen Wissenschaften, Band 196, Springer-Verlag, New
  York-Heidelberg, 1972. \MR{0373040}

\bibitem{DS2008}
A.-L. Basdevant and A.~Singh, \emph{On the speed of a cookie random walk},
  Probab. Theory Related Fields \textbf{141} (2008), no.~3-4, 625--645.
  \MR{2391167}

\bibitem{BW2003}
I.~Benjamini and D.~B. Wilson, \emph{Excited random walk}, Electron. Comm.
  Probab. \textbf{8} (2003), 86--92. \MR{1987097}

\bibitem{CHN2021}
A.~Collevecchio, K.~Hamza, and T.-M. Nguyen, \emph{Long range one-cookie random
  walk with positive speed}, Stochastic Process. Appl. \textbf{142} (2021),
  no.~12, 462--478.

\bibitem{DP2017}
B.~Davis and J.~Peterson, \emph{Excited random walks with non-nearest neighbor
  steps}, J. Theoret. Probab. \textbf{30} (2017), no.~4, 1255--1284.
  \MR{3736173}

\bibitem{Holander}
F.~den Hollander, \emph{Probability theory: The coupling method}, Leiden
  University, 2012.

\bibitem{DK2012}
D.~Dolgopyat and E.~Kosygina, \emph{Scaling limits of recurrent excited random
  walks on integers}, Electron. Commun. Probab. \textbf{17} (2012), no. 35, 14.
  \MR{2965748}

\bibitem{F1971}
W.~Feller, \emph{An introduction to probability theory and its applications.
  {V}ol. {II}}, second ed., John Wiley \& Sons, Inc., New York-London-Sydney,
  1971. \MR{0270403}

\bibitem{Harville}
D.~A. Harville, \emph{Matrix algebra from a statistician's perspective},
  Springer-Verlag, New York, 1997. \MR{1467237}

\bibitem{HW2016}
W.~Hong and H.~Wang, \emph{Branching structures within random walks and their
  applications}, Branching processes and their applications, Lect. Notes Stat.,
  vol. 219, Springer, 2016, pp.~57--73. \MR{3587981}

\bibitem{JS1967}
A.~Joffe and F.~Spitzer, \emph{On multitype branching processes with {$\rho
  \leq 1$}}, J. Math. Anal. Appl. \textbf{19} (1967), 409--430. \MR{212895}

\bibitem{K1991}
S.~V. Kaverin, \emph{Refinement of limit theorems for critical branching
  processes with emigration}, Teor. Veroyatnost. i Primenen. \textbf{35}
  (1990), no.~3, 570--575. \MR{1091216}

\bibitem{Kesten}
H.~Kesten, M.~V. Kozlov, and F.~Spitzer, \emph{A limit law for random walk in a
  random environment}, Compositio Math. \textbf{30} (1975), 145--168.
  \MR{380998}

\bibitem{KM2011}
E.~Kosygina and T.~Mountford, \emph{Limit laws of transient excited random
  walks on integers}, Ann. Inst. Henri Poincar\'{e} Probab. Stat. \textbf{47}
  (2011), no.~2, 575--600. \MR{2814424}

\bibitem{KMP2022}
E.~Kosygina, T.~Mountford, and J.~Peterson, \emph{Convergence of random walks
  with {M}arkovian cookie stacks to {B}rownian motion perturbed at extrema},
  Probab. Theory Related Fields \textbf{182} (2022), no.~1-2, 189--275.
  \MR{4367948}

\bibitem{KP2017}
E.~Kosygina and J.~Peterson, \emph{Excited random walks with {M}arkovian cookie
  stacks}, Ann. Inst. Henri Poincar\'{e} Probab. Stat. \textbf{53} (2017),
  no.~3, 1458--1497. \MR{3689974}

\bibitem{KZ2013}
E.~Kosygina and M.~P.~W. Zerner, \emph{Excited random walks: results, methods,
  open problems}, Bull. Inst. Math. Acad. Sin. (N.S.) \textbf{8} (2013), no.~1,
  105--157. \MR{3097419}

\bibitem{Kosygina2014}
\bysame, \emph{Excursions of excited random walks on integers}, Electron. J.
  Probab. \textbf{19} (2014), no. 25, 25. \MR{3174837}

\bibitem{Lindvall}
T.~Lindvall, \emph{Lectures on the coupling method}, Dover Publications, Inc.,
  Mineola, NY, 2002. \MR{1924231}

\bibitem{Menshikov2012}
M.~Menshikov, S.~Popov, A.~F. Ram\'{\i}rez, and M.~Vachkovskaia, \emph{On a
  general many-dimensional excited random walk}, Ann. Probab. \textbf{40}
  (2012), no.~5, 2106--2130. \MR{3025712}

\bibitem{Norris}
J.~R. Norris, \emph{Markov chains}, Cambridge Series in Statistical and
  Probabilistic Mathematics, vol.~2, Cambridge University Press, Cambridge,
  1998. \MR{1600720}

\bibitem{Hofstad2010}
R.~van~der Hofstad and M.~Holmes, \emph{Monotonicity for excited random walk in
  high dimensions}, Probab. Theory Related Fields \textbf{147} (2010), no.~1-2,
  333--348. \MR{2594356}

\bibitem{V1977}
V.~A. Vatutin, \emph{A critical {G}alton-{W}atson branching process with
  immigration}, Teor. Verojatnost. i Primenen. \textbf{22} (1977), no.~3,
  482--497. \MR{0461694}

\bibitem{VZ1993}
V.~A. Vatutin and A.~M. Zubkov, \emph{Branching processes. {II}}, J. Soviet
  Math. \textbf{67} (1993), no.~6, 3407--3485. \MR{1260986}

\bibitem{V1987}
G.~V. Vinokurov, \emph{On a critical {G}alton-{W}atson branching process with
  emigration}, Teor. Veroyatnost. i Primenen. \textbf{32} (1987), no.~2,
  378--382. \MR{902769}

\bibitem{Z2005}
M.~P.~W. Zerner, \emph{Multi-excited random walks on integers}, Probab. Theory
  Related Fields \textbf{133} (2005), no.~1, 98--122. \MR{2197139}

\end{thebibliography}
\end{document}